\documentclass[11pt]{article}
\usepackage{amsmath}
\usepackage{amssymb,amsmath}
\usepackage{mathrsfs}
\usepackage{graphicx}
\usepackage{color}
\usepackage{esint}%
\usepackage[linktoc=page]{hyperref}
\usepackage[numbers,sort&compress]{natbib}
\usepackage{tocloft}

\def\dive{\mathop{\rm div}\nolimits}

\def\C{\mathop{\bf C\kern 0pt}\nolimits}
\def\DD{\mathop{\bf D\kern 0pt}\nolimits}
\def\K{\mathop{\bf K\kern 0pt}\nolimits}
\def\N{\mathop{\bf N\kern 0pt}\nolimits}
\def\Q{\mathop{\bf Q\kern 0pt}\nolimits}
\def\R{\mathop{\bf R\kern 0pt}\nolimits}

\newcommand{\beq}{\begin{equation}}
\newcommand{\eeq}{\end{equation}}
\newcommand{\ben}{\begin{eqnarray}}
\newcommand{\een}{\end{eqnarray}}
\newcommand{\beno}{\begin{eqnarray*}}
\newcommand{\eeno}{\end{eqnarray*}}

\newtheorem{Theorem}{Theorem}[section]
\newtheorem{Definition}[Theorem]{Definition}

\newtheorem{Lemma}[Theorem]{Lemma}
\newtheorem{Corollary}[Theorem]{Corollary}
\newtheorem{Remark}[Theorem]{Remark}

\numberwithin{equation}{section}

\allowdisplaybreaks \numberwithin{equation} {section}

\begin{document}
\title{NON-UNIQUENESS OF WEAK SOLUTIONS TO 2D GENERALIZED NAVIER-STOKES EQUATIONS  }
      \author{
       \\[2mm]
{\small  XINLIANG LI  \ \ \ \ \  \ ZHONG TAN}}
         \date{}
         \maketitle
\noindent{\bf Abstract.} We study the non-uniqueness of weak solutions for the two-dimensional hyper-dissipative Navier-Stokes equations in the super-critical spaces $L_{t}^{\gamma}L_{x}^{p}$ when $\alpha\in[1,\frac{3}{2})$, and obtain the conclusion that the non-uniqueness of the weak solutions at the endpoint $(\gamma,p)=(\infty, \frac{2}{2\alpha-1})$ is sharp in view of the generalized Lady\v{z}enskaja-Prodi-Serrin condition by using a different  spatial-temporal building block from [Cheskidov-Luo, Ann. PDE, 9:13 (2023)] and taking advantage of the intermittency of the temporal concentrated function $g_{(k)}$ in an almost optimal way. Our results recover the above 2D non-uniqueness conclusion and extend to the hyper-dissipative case $\alpha \in(1,\frac{3}{2})$.\ \ \

\vskip   0.2cm \footnotetext[1]{College of Mathematical Sciences, Shenzhen University, Shenzhen 518061, P. R. China; College of Physics and
Optoelectronic Engineering, Shenzhen University, Shenzhen 518061, P. R. China. Email address, Xinliang Li: lixinliangmaths@163.com.}

\vskip   0.2cm \footnotetext[2]{College of Mathematical Sciences, Xiamen University, Xiamen 361005, P. R. China;
Shenzhen Research Institute of Xiamen University, Shenzhen 518057, P. R. China.
Email address, Zhong Tan: tan85@xmu.edu.cn.}

	\setlength{\baselineskip}{20pt}
\vspace{-5mm}
\begin{center}
\tableofcontents
\end{center}

\begin{center}
\section{Introduction And Main Results}\label{A}
\end{center}
\subsection{Inroduction.}\ \ \ \
In this paper, we consider the following 2D generalized Navier-Stokes equations (gNSE):
\[
\left\{\begin{array}{l}
\partial_{t} u+\nu(-\Delta)^{\alpha} u+(u \cdot \nabla) u+\nabla\mathsf{P}=0 , \label{1.1}\tag{1.1}\\
\operatorname{div} u=0,\\
u|_{t=0}=u_{0},
\end{array}\right.
\]
on the torus $\mathbb{T}^{2}:=\mathbb{R}^{2}/\mathbb{Z}^{2}$, where $u:\mathbb{T}^{2}\times[0,T] \rightarrow\mathbb{R}^{2}$ is the unknown velocity and $\mathsf{P}: \mathbb{T}^{2}\times[0,T]\rightarrow \mathbb{R}$ is the scalar pressure, $\nu>0$ is the constant viscosity coefficient, $\alpha$ is called viscosity exponent in the literature, and $(-\Delta)^{\alpha}$ is the fractional Laplacian defined by the Fourier transform on the flat torus as
$$
\mathcal{F}\left((-\Delta)^{\alpha} u\right)(\xi)=|\xi|^{2 \alpha} \mathcal{F}(u)(\xi), \quad \xi \in \mathbb{Z}^{2}.
$$
When $\alpha < 1$, \eqref{1.1} is known as hypo-dissipative Navier-Stokes equations, and  $\alpha > 1$ is hyper-dissipative Navier-Stokes equations. When $\alpha = 1$, \eqref{1.1} turns to be the 2D classical Navier-Stokes equations (NSE), for which the Leray-Hopf solution is the unique smooth solution in the physical sense. J.-L. Lions \cite{JL69} first  studied  the system \eqref{1.1} and proved the Leray-Hopf solution is unique to the $d\geq2$ dimensional gNSE with any $\alpha \geq \frac{1}{2}+\frac{d}{4}$, which also guarantees global regularity with smooth initial data presented in Wu \cite{JW03}.

For the 2D hyper-dissipative NSE \eqref{1.1}, we mainly consider the mixed Sobolev space, and the space $\mathbb{X}=L_{t}^{\gamma} L_{x}^{p}$ is called critical if its norm $\|\cdot\|_{\mathbb{X}}$ is invariant under the natural scaling:
\begin{equation*}
u(t, x) \mapsto \lambda^{2 \alpha-1} u\left(\lambda^{2 \alpha} t, \lambda x\right), \quad P(t, x) \mapsto \lambda^{4 \alpha-2} P\left(\lambda^{2 \alpha} t, \lambda x\right), \label{1.4}\tag{1.2}
\end{equation*}
where the exponents $(\gamma, p)$ satisfy
\begin{equation*}
\frac{2 \alpha}{\gamma}+\frac{2}{p}=2 \alpha-1 \label{1.5}\tag{1.3}
\end{equation*}
known as generalized Lady\v{z}enskaja-Prodi-Serrin (gLPS) condition.

In particular, the mixed Lebesgue space $L_{t}^{\gamma} L_{x}^{p}$ is critical for the classical NSE ($\alpha=1$), when the exponents $(\gamma, p)$ satisfy the well-known Lady\v{z}enskaja-Prodi-Serrin (LPS) condition $\frac{2}{\gamma}+\frac{2}{p}=1$(see \cite{PRO59,SER62,LAD67}, for the case $d=2$). When $\frac{2}{\gamma}+\frac{2}{q}<1$,  the space $L_{t}^{\gamma} L_{x}^{p}$ is called sub-critical and super-critical when $\frac{2}{\gamma}+\frac{2}{q}>1$.

The known uniqueness results can be concisely summarized as stated in Theorem 1.3 of \cite{AC22}. Specifically, the findings reported in \cite{EB72,GF00,PL01} imply that any weak solution to the 2D NSE in the distributional sense, as defined in Definition \ref{A.1} below, within the (sub)critical spaces $L_{t}^{\gamma} L_{x}^{p}$ with the condition $\frac{2}{\gamma}+ \frac{2}{p} \leq 1$ for some $\gamma, p \in [1, \infty]$, is automatically the unique regular Leray-Hopf solution on $\mathbb{T}^{2}$. Furthermore, there also exist many uniqueness results for the 2D hyper-dissipative NSE $(\alpha>1)$, either under the gLPS condition \eqref{1.5} or  subcritical spaces, such as \cite{ST98} for the one endpoint space $C_{t} L_{x}^{\frac{2}{2\alpha-1}}$, and \cite{ST03} for the Besov space $\dot{B}_{p, q}^{\frac{2}{p}+1-2 \alpha}$.

In the super-critical space, there remains many non-uniqueness questions unsolved. A recent noteworthy contribution by Cheskidov and Luo \cite{AC22} has established the sharp non-uniqueness of the 2D NSE in proximity to the critical endpoint $(\gamma, p) = (2, \infty)$ within the LPS condition. This proof heavily relies on the exploitation of temporal intermittency in the Convex integration framework. Additional references \cite{AC21,ACX22} demonstrate the application of temporal intermittency to transport equations, while \cite{YLZ22,YLZZ22} explores its implications for 3D Magneto-Hydrodynamic (MHD) equations. Separately, Li, Qu, Zeng and Zhang \cite{YLP22} have presented non-unique $L_{t}^{2}C_{x}$ weak solutions for the hypo-viscous NSE $(\alpha<1)$ in all dimensions $d\geq2$ which implies that $\alpha = 1$ is the sharp viscosity threshold for the well-posedness in the space $L_{t}^{2}C_{x}$. Regarding the other endpoint $(\gamma, p) = (\infty, 2)$, Luo and Qu \cite{LTP20} have established that $C_{t}L_{x}^{2}$ weak solutions lack uniqueness when $\alpha<1$ in 2D case. Recently, Cheskidov and Luo \cite{ACL22} proved that it is not feasible to construct non-unique $C_{t}L_{x}^{2}$ weak solutions for the 2D NSE, indicating that $\alpha = 1$ serves as the sharp viscosity threshold for well-posedness within the space $C_{t}L_{x}^{2}$.

In summary, the current research indicates that there exist non-unique results in the super-critical space for 2D generalized Navier-Stokes equations $\alpha\leq1$. However, the non-unique results for  $\alpha>1$ in 2D case remain unknown. In this study, we give the first non-uniqueness result in the Theorem \ref{A.2} below when $\alpha\in(1,\frac{3}{2})$.

\subsection{Spatial-temporal Convex Integration Technique.}\ \ \ \
The construction of non-unique weak solutions in the present work relies on Convex integration techniques, which was introduced to fluid dynamics by De Lellis and  Sz\'{e}kelyhidi Jr. in \cite{CD09}. Since then, it was developed across a series of works \cite{ TB15, TB16,CD09, CD13, CD14} ultimately leading to the resolution of the flexible part of Onsager's conjecture for the 3D Euler equations by Isett \cite{PI18} and Buckmaster-De Lellis-Sz\'{e}kelyhidi-Vicol \cite{TB19}. When considering the viscous fluids, Buckmaster and Vicol \cite{BV19} constructed non-unique weak solutions $C_{t}L_{x}^{2}$ of the 3D NSE by using the method of intermittent Convex integration. Following this significant advancement, there have been many subsequent results in viscous settings including extensions of \cite{BV19} to the hyper-dissipative case \cite{LT20, YL22}, partial regularity in time \cite{BCV22}, high dimensions $d \geq 4$ \cite{Luo19}, MHD equations \cite{YLZ22,YLZZ22} and Hall-MHD equations \cite{Dai21} and so on. For the non-uniqueness of Leray-Hopf solutions to 3D hypo-dissipative NSE please refer to \cite{CDLDR18, DR19}.

In 2D case, due to the principle that any two non-parallel lines on the same plane must intersect, how to set up different building blocks in the Convex integration framework without interfering or intersecting with each other is one of the key considerations. Additionally, in order to achieve hyper-viscosity $(\alpha>1)$, using only the spatial intermittent Convex integration method \cite{LTP20} is not sufficient to deal with the difficulties brought about by the hyper-dissipative error. Therefore, it is necessary to introduce an additional temporal concentration function $g_{(k)}$ in two-dimensional space, for which the validity has been verified in \cite{AC21,AC22,ACL22,ACX22}. This study also requires the use of space-time intermittent Convex integration method, and will focus on developing the crucial role of temporal concentration functions $g_{(k)}$ so that even if the different spatial building blocks $\mathbf{W}_{k}$ intersects with each other, the newly generated ``Building blocks'' $g_{(k)}\mathbf{W}_{k}$ will not intersect under the effect of non-intersecting temporal building blocks by temporal shifts $t_{k}$. We here vividly refer to the temporal concentration function $g_{(k)}$ as a ``switch function" and we will utilize temporal intermittency in an almost optimal way to achieve the sharp non-uniqueness at the endpoint cases of gLPS condition \eqref{*}, namely $(\infty, \frac{2}{2\alpha-1})$. As will be demonstrated subsequently, in the case of extreme dissipativity where $\alpha$ approximates $\frac{3}{2}$, the appropriate temporal intermittency is approximately equivalent to 2D spatial intermittency.

{\bf Notations.} The notation $a \lesssim b$ means that $a \leq C b$ for some non-negative constant $C$. The mean of $u \in L^1\left(\mathbb{T}^2\right)$ is given by $\fint_{\mathbb{T}^2} u d x=\left|\mathbb{T}^2\right|^{-1} \int_{\mathbb{T}^2} u d x$, where $|\cdot|$ denotes the Lebesgue measure. For simplicity, for $p, \gamma \in[1, \infty]$, we denote
$$
L_{t}^{\gamma}:=L^{\gamma}(0, T), \quad L_{x}^{p}:=L^{p}\left(\mathbb{T}^{2}\right),
$$
where $L_{t}^{\gamma} L_{x}^{p}$ denotes the usual Banach space $L^{\gamma}\left(0, T ; L^{p}\left(\mathbb{T}^{2}\right)\right)$. In particular, we write $L_{t, x}^{p}:=L_{t}^{p} L_{x}^{p}$ for brevity. Let
$$
\|u\|_{W_{t, x}^{N, p}}:=\sum_{0 \leq m+|\xi| \leq N}\left\|\partial_{t}^{m} \nabla^{\xi} u\right\|_{L_{t, x}^{p}}, \quad\|u\|_{C_{t, x}^{N}}:=\sum_{0 \leq m+|\xi| \leq N}\left\|\partial_{t}^{m} \nabla^{\xi} u\right\|_{C_{t, x}}
$$
where $\xi=\left(\xi_{1}, \xi_{2} \right)$ is the multi-index and $\nabla^{\xi}:=\partial_{x_{1}}^{\xi_{1}} \partial_{x_{2}}^{\xi_{2}}$. For any Banach space $X$, we use $C([0, T]; X)$ to denote the space of continuous functions from $[0, T]$ to $X$ with norm $\|u\|_{C_{t} X}:=\sup _{t \in[0, T]}\|u(t)\|_{X}$,
and $B_{p, q}^{s}\left(\mathbb{T}^{2}\right)$ to denote the Besov space with the norm
$$
\|u\|_{B_{p, q}^{s}\left(\mathbb{T}^{2}\right)}=\left(\sum_{j \geq -1}\left|2^{j s}\left\|\Delta_{j} u\right\|_{L^{p}\left(\mathbb{T}^{2}\right)}\right|^{q}\right)^{\frac{1}{q}}
$$
where $\left\{\Delta_{j}\right\}_{j \in \mathbb{Z}}$ is the inhomogeneous frequency localization operators.

For any $\mathscr{A} \subseteq[0, T]$, we introduce the $\varepsilon_{*}$-neighborhood of $\mathscr{A}$ as
$$
\mathscr{N}_{\varepsilon_{*}}(\mathscr{A}):=\left\{t \in[0, T]: \exists\ s \in \mathscr{A} \text {, s.t. }|t-s| \leq \varepsilon_{*}\right\},
$$
where $\varepsilon_{*}>0$. Additionally, we set the viscosity coefficient $\nu=1$ for brevity.

\subsection{Main Results.}\ \ \ \
Before presenting the main results, we introduce the concept of weak solutions in the distributional sense for the equations \eqref{1.1}.

\begin{Definition} \label{A.1}
(Weak solution) For $\alpha\geq1$, given any weakly divergence-free datum $u_{0} \in L^{2}\left(\mathbb{T}^{2}\right)$, $u \in$ $L^{2}\left([0, T] \times \mathbb{T}^{2}\right)$ is a weak solution for the hyper-dissipative NSE \eqref{1.1}, if $u$ is divergence-free for all $t \in[0, T]$ and satisfies
\begin{align*}
\int_{\mathbb{T}^{2}} u_{0} \Phi(0, \cdot) d x=-\int_{0}^{T} \int_{\mathbb{T}^{2}} u\left(\partial_{t} \Phi-(-\Delta)^{\alpha} \Phi+(u \cdot \nabla)\Phi\right) d x d t ,
\end{align*}
where $\Phi \in C_{0}^{\infty}\left([0, T) \times \mathbb{T}^{2}\right)$ is the divergence-free test function.
\end{Definition}

The main result of this paper is formulated in Theorem \ref{A.2} below.
\begin{Theorem} \label{A.2}
For $\alpha \in[1,\frac{3}{2})$. If one has a smooth, divergence-free and mean-free vector field $\tilde{u}=\tilde{u}(t,x)$ on $[0, T] \times \mathbb{T}^{2}$, and $(\gamma, p)$ satisfying
\begin{align*}
\frac{4\alpha-4}{\gamma}+\frac{2}{p}>2\alpha-1. \label{*}\tag{1.4}
\end{align*}
Then, there exists $\beta^{\prime} \in(0,1)$, such that for any given $\varepsilon_*>0$, there exists a velocity field $u$ such that the following holds:

(i) Weak solution: $u$ is a weak solution to \eqref{1.1} in the sense of Definition \ref{A.1} with zero spatial mean.

(ii) Regularity: $u\in H_{t, x}^{\beta^{\prime}} \cap L_t^\gamma L_x^{p}$.

(iii) $\varepsilon_*$-neighborhood of the temporal support: $\operatorname{supp}_t u \subseteq \mathscr{N}_{\varepsilon_{*}}\left(\operatorname{supp}_t \tilde{u} \right)$.

(iv) $\varepsilon_*$-close between the solution $u$ and the given $\widetilde{u}$ in $L_t^1 L_x^2\cap L_t^\gamma L_x^p$,
$$\|u-\tilde{u}\|_{L_t^1 L_x^2\cap L_t^\gamma L_x^p} \leq \varepsilon_*.$$
\end{Theorem}

The following non-uniqueness of weak solutions to \eqref{1.1} is actually a direct consequence of Theorem \ref{A.2} when $\alpha\in[1, \frac{3}{2})$. Please refer to \cite{AC22,YL22} for the proof of relevant details.
\begin{Corollary}\label{A.3}
(Non-uniqueness for 2D hyper-dissipative NSE) Let $\alpha \in[1,\frac{3}{2})$. For any weak solution $\widetilde{u}$ to \eqref{1.1}, there exists a different weak solution $u \in L_{t}^{\gamma} L_{x}^{p}$ to \eqref{1.1} with the same initial data, where $(\gamma, p)$ satisfies \eqref{*}.

Furthermore, for any divergence-free initial data belonging to $L_{x}^{2}$, there exist infinitely many weak solutions in $L_{t}^{\gamma} L_{x}^{p}$ to \eqref{1.1}.
\end{Corollary}

\subsection{Comments On Main Results.}\ \ \ \
In this subsection, we will share some comments on the the key observations drawn from our results.

{\bf (1) Strong non-uniqueness for the high viscosity $\alpha\geq1$.} According to the classical result \cite{JL69}, there exists a unique Leray-Hopf weak solution to \eqref{1.1} when $\alpha\geq1$, it seems a novel finding to get the non-unique result when $\alpha=1$ or even greater than 1 in Theorem \ref{A.2}. In order to explore the non-uniqueness of weak solutions in the super-critical space \eqref{*} when $\alpha\geq1$, it is necessary to fully exploit the temporal intermittency. To this end, we establish  the relationship between the intermittency of the temporal concentrated function $g_{(k)}$, the viscosity exponent $\alpha$ and the regularity of the weak solutions, see \eqref{4.1} and gLPS condition \eqref{1.5} for details. Our results recover the two-dimensional results in \cite{ACL22} obtained by Cheskidov and Luo when $\alpha=1, \gamma=\infty$ in \eqref{*}, and are the first non-uniqueness results when $\alpha>1$ to \eqref{1.1}.

{\bf (2) Sharp non-uniqueness at the endpoint of gLPS condition.} In the remarkable papers \cite{ACL22}, Cheskidov-Luo first proved the sharp non-uniqueness at the  endpoint cases  $(\gamma, p)=(\infty, 2)$ for the 2D NSE.

 In view of the well-posedness result in critical space $C_{t}L_{x}^{2 /(2 \alpha-1)}$ by \cite{ST98} when $\alpha\geq1$, Theorem \ref{A.2} provides the sharp non-uniqueness for the hyper-dissipative NSE \eqref{1.1} at the endpoint $(\gamma,p)=(\infty, \frac{2}{2\alpha-1})$.  Therefore, we extend the sharp non-uniqueness result in \cite{ACL22} to the 2D hyper-dissipative NSE for $\alpha \in(1,\frac{3}{2})$. In particularly, the viscosity exponent $\alpha=\frac{3}{2}$ seems the upper limit in the super-critical spaces due to $\frac{2}{2\alpha-1}>1$.

{\bf (3) Sharp non-uniqueness at the other endpoint of gLPS condition.} For the other endpoint case $(\gamma,p)=(\frac{2\alpha}{2\alpha-1}, \infty)$, the detailed proof will be presented in our future research work \cite{DL24}. At present, due to the $L_{t, x}^{2}$-criticality of space-time Convex integration method, whether the non-uniqueness of the weak solutions is sharp for non-endpoint cases in view of the gLPS condition \eqref{1.5} is unknown,
we are looking forward to better technical means for processing.

\begin{center}
\section{Main Iteration And Mollification procedure}\label{B}
\end{center}\ \ \ \
In this section, we first present the main iteration of the velocity and Reynolds stress, which is the heart of the proof of our main theorem.\\
\subsection{Main Iteration}\ \ \ \ We consider the approximate solutions to the following 2D Navier-Stokes-Reynolds system for each  integer $q \geq 0$,
\[
\left\{\begin{array}{l}\label{2.1}
\partial_{t} u_{q}+(-\Delta)^{\alpha} u_{q}+\operatorname{div}\left(u_{q} \otimes u_{q}\right)+\nabla\mathsf{P}_{q}=\operatorname{div} \mathring{\mathsf{R}}_{q},  \tag{2.1}\\
\operatorname{div} u_{q}=0,
\end{array}\right.
\]
where $\mathring{\mathsf{R}}_{q}: \mathbb{T}^{2}\times[0,T]\rightarrow \mathcal{S}_{0}^{2\times2}$ is a $2 \times 2$ symmetric traceless  matrix, known as Reynolds stress in the literature, where the pressure term can be recoverd from the elliptic equation by taking the divergence of \eqref{2.1}. Specifically,
$$
\Delta \mathsf{P}_{q}=\operatorname{div} \operatorname{div}(\mathring{\mathsf{R}}_{q}-u_{q} \otimes u_{q}),
$$
which, in conjunction with the standard zero spatial mean condition $\fint_{\mathbb{T}^2} \mathsf{P} dx= 0$, uniquely determines the pressure.

To measure the size of the relaxation solutions $\left(u_{q}, \mathring{\mathsf{R}}_{q}\right)$ for $q \in \mathbb{N}$, we employ two crucial parameters: the frequency parameter $\lambda_{q}$ and the amplitude parameter $\delta_{q+2}$:
\begin{equation*}
\lambda_{q}=a^{\left(b^{q}\right)}, \quad \delta_{q+2}=\lambda_{q+2}^{-2 \beta}. \label{2.2}\tag{2.2}
\end{equation*}
where $a\in\mathbb{N}$ is a large integer, the parameter $\beta>0$ represents regularity and $b\in 2 \mathbb{N}$ satisfies
\begin{equation*}
b>\frac{1000}{\varepsilon }, \quad 0<\beta<\frac{1}{100 b^{2}}, \label{2.3}\tag{2.3}
\end{equation*}
where $\varepsilon \in \mathbb{Q}_{+}$ is sufficiently small such that
\begin{equation*}
\varepsilon \leq \frac{1}{20} \min \left\{3-2\alpha, \frac{4 \alpha-4}{\gamma}+\frac{2}{p}-(2 \alpha-1)\right\} \quad \text { and } \quad b \varepsilon \in \mathbb{N}. \label{2.4}\tag{2.4}
\end{equation*}

The purpose of the iteration is to prove that when $q$ approaches infinity, the Reynolds stress disappears in an appropriate space and the limit value of $u_{q}$ is the solution of the initial equations \eqref{1.1}.  This approach is formally outlined in the subsequent iterative steps:
\begin{align*}
& \left\|u_{q}\right\|_{C_{t,x}^{1}} \lesssim \lambda_{q}^{7},  \label{2.6}\tag{2.5}\\
& \left\|\mathring{\mathsf{R}}_{q}\right\|_{L_{t, x}^{1}} \leq \delta_{q+1},\label{2.10} \tag{2.6}\\
& \left\|\mathring{\mathsf{R}}_{q}\right\|_{C_{t,x}^{1}} \lesssim \lambda_{q}^{16},  \label{2.8}\tag{2.7}
\end{align*}
where the implicit constants $C$ are independent of $q$.

Our main theorem relies on the following
crucial iteration lemma formulated below.
\begin{Lemma}\label{B.3}(Iteration lemma). Let $\alpha \in[1,\frac{3}{2})$ and $(p, \gamma)$ satisfy \eqref{*}.
Suppose that $\left(u_{q}, \mathring{\mathsf{R}}_{q}\right)$ is a smooth solution to \eqref{2.1} satisfying conditions \eqref{2.6}-\eqref{2.8}. Then, there exists another solution $\left(u_{q+1}, \mathring{\mathsf{R}}_{q+1}\right)$ to \eqref{2.1} fulfills conditions \eqref{2.6}-\eqref{2.8} with $q+1$ replacing $q$. Furthermore, we have
\begin{align*}
& \left\|u_{q+1}-u_{q}\right\|_{L_{t, x}^{2}} \leq \delta_{q+1}^{\frac{1}{2}},  \label{2.11}\tag{2.8}\\
& \left\|u_{q+1}-u_{q}\right\|_{L_{t}^{1} L_{x}^{2}} \leq \delta_{q+2}^{\frac{1}{2}},  \label{2.12}\tag{2.9}\\
& \left\|u_{q+1}-u_{q}\right\|_{L_{t}^{\gamma} L_{x}^{p}} \leq \delta_{q+2}^{\frac{1}{2}}, \label{2.13}\tag{2.10}
\end{align*}
and
\begin{equation*}
\operatorname{supp}_{t}\left(u_{q+1}, \mathring{\mathsf{R}}_{q+1}\right) \subseteq \mathscr{N}_{\delta_{q+2}^{\frac{1}{2}}}\left(\operatorname{supp}_{t}\left(u_{q}, \mathring{\mathsf{R}}_{q}\right)\right). \label{2.14}\tag{2.11}
\end{equation*}
\end{Lemma}

\subsection{Mollification Procedure}\ \ \ \
Let $\Phi_\epsilon$ and $\widetilde{\Phi}_\epsilon$ be standard mollifiers on $\mathbb{R}^2$ and $\mathbb{R}$, respectively, and $\operatorname{supp}_{t} \widetilde{\Phi}_{\epsilon} \subseteq$ $(-\epsilon, \epsilon)$. Then, we can mollify $u_{q}$ and $\mathring{\mathsf{R}}_{q}$ given in Lemma \ref{B.3} as
\begin{align*}
u_{\ell}  =\left(u_{q} *_{x} \Phi_{\ell}\right) *_{t} \widetilde{\Phi}_{\ell},  \quad
\mathring{\mathsf{R}}_{\ell} =\left(\mathring{\mathsf{R}}_{q} *_{x} \Phi_{\ell}\right) *_{t} \widetilde{\Phi}_{\ell} \label{212}\tag{2.12}
\end{align*}
where the scale of mollification is given by
\begin{equation*}
\ell=\lambda_{q}^{-20}.\label{213} \tag{2.13}
\end{equation*}

Since $\left(u_{q}, \mathsf{P}_{q}, \mathring{\mathsf{R}}_{q}\right)$ solves \eqref{2.1}, we know that $\left(u_{\ell}, \mathsf{P}_{\ell}, \mathring{\mathsf{R}}_{\ell}\right)$ solves
\[
\left\{\begin{array}{l}
\partial_{t} u_{\ell}+\dive\left(u_{\ell} \otimes u_{\ell}\right)+\nabla \mathsf{P}_{\ell}+(-\Delta)^{\alpha} u_{\ell}=\dive\left(\mathring{\mathsf{R}}_{\ell}+\mathring{\mathsf{R}}_{\text {com}}\right),  \label{214}\tag{2.14}\\
\dive u_{\ell}=0,
\end{array}\right.
\]
where we can choose
\begin{align*}
&\mathsf{P}_{\ell}=\left(\mathsf{P}_{q} *_{x} \Phi_{\ell}\right) *_{t} \widetilde{\Phi}_{\ell}+\left|u_{\ell}\right|^{2}-\left(\left|u_{q}\right|^{2} *_{x} \Phi_{\ell}\right) *_{t} \widetilde{\Phi}_{\ell},  \label{215}\tag{2.15}\\
&\mathring{\mathsf{R}}_{\text {com }} =\left(u_{\ell} \mathring{\otimes} u_{\ell}\right)-\left(\left(v_{q} \mathring{\otimes} u_{q}\right) *_{x} \Phi_{\ell}\right) *_{t} \widetilde{\Phi}_{\ell}. \label{216}\tag{2.16}
\end{align*}

Using the inductive assumptions \eqref{2.6}-\eqref{2.8}, we have
\begin{gather*}
\left\|u_{\ell}\right\|_{C_{t, x}^{N}} \lesssim \lambda_{q}^{7} \ell^{-N+1} \lesssim \ell^{-N}, \label{217}\tag{2.17}\\
\left\|\mathring{\mathsf{R}}_{\ell}\right\|_{C_{t, x}^{N}} \lesssim \lambda_{q}^{16} \ell^{-N+1} \lesssim \ell^{-N}, \label{218}\tag{2.18}\\
\left\|\mathring{\mathsf{R}}_{\ell}\right\|_{L_{t,x}^{1}} \leq\left\|\mathring{\mathsf{R}}_{q}\right\|_{L_{t,x}^{1}} \leq  \delta_{q+1},  \label{219}\tag{2.19}\\
\left\|u_{\ell}-u_{q}\right\|_{L_{t}^{\gamma} L_{x}^{p}} \lesssim\left\|u_{\ell}-u_{q}\right\|_{L_{t,x}^{\infty}} \lesssim \ell\left\|u_{q}\right\|_{C_{t, x}^{1}} \lesssim \lambda_{q}^{-13} . \label{220}\tag{2.20}
\end{gather*}

Moreover,
$$
\begin{aligned}
\left\|\mathring{\mathsf{R}}_{\text {com}}\right\|_{L_{t,x}^{\infty}} & \lesssim \ell \|u_{\ell} \otimes u_{\ell} \|_{C_{t, x}^{1}} \lesssim \ell \lambda_{q}^{14}, \\
\left\|\mathring{\mathsf{R}}_{\text {com}}\right\|_{C_{t, x}^{N}} & \lesssim \ell^{-N+1}\left\|u_{\ell} \mathring{\otimes} u_{\ell}\right\|_{C_{t, x}^{1}} \lesssim \ell^{-N+1} \lambda_{q}^{14}.
\end{aligned}
$$
Thus, for
$$
\mathring{\mathsf{R}}_{\ell}^{*} \stackrel{\text { def. }}{=} \mathring{\mathsf{R}}_{\ell}+\mathring{\mathsf{R}}_{\text {com}},
$$
we have
\begin{gather*}
\left\|\mathring{\mathsf{R}}_{\ell}^{*}\right\|_{L_{t,x}^{1} } \leq  \delta_{q+1}+\ell \lambda_{q}^{14} \leq  \delta_{q+1},  \label{221}\tag{2.21}\\
\left\|\mathring{\mathsf{R}}_{\ell}^{*}\right\|_{C_{t, x}^{N}} \lesssim \ell^{-N}+\ell^{-N+1} \lambda_{q}^{10} \lesssim \ell^{-N}, \label{222}\tag{2.22}
\end{gather*}
where we use the fact that by \eqref{2.3} and \eqref{213}, it holds
$$
\ell \lambda_{q}^{14} \leq \delta_{q+1}.
$$

{\centering
\section{Constructions Of Velocity Perturbations\label{C}}}\ \
In this section,  we are primarily focus on constructing appropriate velocity perturbations and selecting appropriate intermittent or oscillatory parameters, to make the corresponding inductive estimates hold in Lemma \ref{B.3}.

In order to approach one endpoint $(\gamma, p)=(\infty, \frac{2}{2\alpha-1})$ with $\alpha \in[1,\frac{3}{2})$, we select the 2D intermittent jets $\mathbf{W}_{k}$ as the primary spatial building blocks (see \eqref{4.5} below). A key characteristic of the intermittent jet is its near-2D intermittency, i.e.,
\begin{equation*}
\left\|\mathbf{W}_{k}\right\|_{L_{t}^{\infty} L_{x}^{1}} \lesssim \lambda^{-1+}. \label{2.19}
\end{equation*}
To effectively manage the hyper-dissipativity error with $(-\Delta)^{\alpha}$ when $\alpha>1$,  it is necessary to introduce the temporal concentration function $g_{(k)}(t)$. Specifically, the appropriate temporal intermittency closely corresponds to $(4 \alpha-4)$-dimensional spatial intermittency, when $\alpha$ approaches $\frac{3}{2}$, the temporal intermittency nearly attains 2D spatial intermittency. Furthermore, we demonstrate the existence of six admissible parameters: $\left(r_{\perp}, r_{\|}, \lambda, \mu, \tau, \sigma \right)$ and provide a precise selection as follows:
\begin{equation*}
r_{\perp}:=\lambda_{q+1}^{-1+2 \varepsilon}, \ r_{\|}:=\lambda_{q+1}^{-1+10 \varepsilon},\ \lambda:=\lambda_{q+1},\ \mu:=\lambda_{q+1}^{2\alpha-1 +4 \varepsilon},\ \tau:=\lambda_{q+1}^{4 \alpha-4+16 \varepsilon}, \ \sigma:=\lambda_{q+1}^{2 \varepsilon} \label{4.1}\tag{3.1}
\end{equation*}
where $\varepsilon$ denotes a sufficiently small constant that fulfills the condition stated in \eqref{2.4}. \\

\subsection{Spatial Building Blocks.}\label{D.3.1}\ \ \ \
Define $\varphi, \psi:\mathbb{R}\rightarrow \mathbb{R}$ as smooth, mean-free functions supported on a ball of radius 1, satisfying the following conditions
\begin{equation*}
\frac{1}{2 \pi} \int_{\mathbb{R}} \varphi^{2}(x) \mathrm{d} x=1\quad \text{and}\quad \frac{1}{2 \pi} \int_{\mathbb{R}} \psi^{2}(x) \mathrm{d} x=1, \quad \operatorname{supp} \varphi,\ \psi \subseteq[-1,1]. \label{4.2}\tag{3.2}
\end{equation*}

The corresponding rescaled cut-off functions are defined as follows:
\begin{equation*}
\varphi_{r_{\|}}(x):=r_{\|}^{-\frac{1}{2}}\varphi\left(\frac{x}{r_{\|}}\right), \quad \psi_{r_{\perp}}(x):=r_{\perp}^{-\frac{1}{2}} \psi\left(\frac{x}{r_{\perp}}\right) .\label{4.3}\tag{3.3}
\end{equation*}
In the scaling process, the functions $\varphi_{r_{\|}}$ and $\psi_{r_{\perp}}$ are supported within the balls of radius $r_{\|}$ and $r_{\perp}$ in $\mathbb{R}$, respectively. By an abuse of notation, we extend the definition of $\varphi_{r_{\|}}$ and $\psi_{r_{\perp}}$ to periodic functions on the torus $\mathbb{T}$.

Let $\Lambda$ be a subset of $\mathbb{S} \cap \mathbb{Q}^{2}$ representing the wavevector set, as specified in Geometric Lemma \ref{G.1}. For each $k \in \Lambda$, let $\left(k, k_{1}\right)$ denotes an orthonormal basis. Then the 2D intermittent jets are defined by
\begin{align*}
\mathbf{W}_{k}:=-\varphi_{r_{\|}}\left(\lambda r_{\perp} N_{\Lambda}(k_{1} \cdot x+\mu t)\right) \psi_{r_{\perp}}^{\prime}\left(\lambda r_{\perp} N_{\Lambda} k\cdot x\right) k_{1}, \label{4.4}\tag{3.4}
\end{align*}
where the quantity $N_{\Lambda}$ is defined according to \eqref{7.2}. By design, $\mathbf{W}_{k}$ is $(\mathbb{T}/\lambda r_{\perp})^2$-periodic in space and $\sigma^{-1}$-periodic in time.

The parameters $r_{\|}$ and $r_{\perp}$ serve as indicators of the concentration effect exerted by intermittent jets. Additionally, $\mu$ represents the temporal oscillation parameter of the flow within the building block.

For brevity, we set
\begin{align*}
& \varphi_{\left(k_{1}\right)}(x):=\varphi_{r_{\|}}\left(\lambda r_{\perp} N_{\Lambda}(k_{1} \cdot x+\mu t)\right)=\varphi_{r_{\|}}\left(\lambda r_{\perp} N_{\Lambda}(x_{k}+\mu t)\right), \\
& \psi_{(k)}(x):=\psi_{r_{\perp}}\left(\lambda r_{\perp} N_{\Lambda} k \cdot x\right)=\psi_{r_{\perp}}\left(\lambda r_{\perp} N_{\Lambda} y_{k}\right),
\end{align*}
where the pair $(x_{k}, y_{k})$ is employed to represent the spatial coordinates in $\mathbb{R}^{2}$. Therefore, the intermittent jets is simplified as
\begin{equation*}
\mathbf{W}_{k}=-\varphi_{\left(k_{1}\right)} \psi_{(k)}^{\prime} k_{1}, \quad k \in \Lambda . \label{4.5}\tag{3.5}
\end{equation*}
Since $\mathbf{W}_{k}$ is not divergence-free, we introduce the corrector $\mathbf{W}_{k}^{c}$  as
\begin{equation*}
\mathbf{W}_{k}^{c}:=\frac{r_{\perp}}{r_{\|}}\varphi_{\left(k_{1}\right)}^{\prime} \psi_{(k)} k. \label{4.6}\tag{3.6}
\end{equation*}
Furthermore, we define periodic potentials $\boldsymbol{\Psi}_{k} \in C_{c}^{\infty}\left(\mathbb{R}^{2}\right)$ as
\begin{equation*}
\boldsymbol{\Psi}_{k} =r_{\perp}\varphi_{\left(k_{1}\right)} \psi_{(k)}. \label{4.7}\tag{3.7}
\end{equation*}
By straightforward computations,
$$
\begin{aligned}
(\lambda r_{\perp}N_{\Lambda})^{-1} \nabla^{\perp} \boldsymbol{\Psi}_{k} & =(\lambda r_{\perp}N_{\Lambda})^{-1}(-\partial_{y_{k}} \boldsymbol{\Psi}_{k} k_{1}+\partial_{x_{k}} \boldsymbol{\Psi}_{k} k) \\
& =\mathbf{W}_{k}+\mathbf{W}_{k}^{(c)}.
\end{aligned}
$$
Also, we have the important identities
\begin{equation*}
\partial_{t}\left|\mathbf{W}_{k}\right|^{2} k_{1}= \mu  \operatorname{div}\left(\mathbf{W}_{k} \otimes \mathbf{W}_{k}\right), \label{4.8}\tag{3.8}
\end{equation*}
and
\begin{equation*}
\partial_{t} \boldsymbol{\Psi}_{k}= \mu \left(k_{1} \cdot \nabla\right) \boldsymbol{\Psi}_{k}. \label{4.9}\tag{3.9}
\end{equation*}

After a direct calculation, we have the following lemma is the crucial estimates pertaining to the 2D intermittent jets.

\begin{Lemma}\label{D.1}(Estimates of 2D intermittent jets). For any $p \in[1, \infty]$, $N, M \in \mathbb{N}$, the following bounds hold,
\begin{align*}
& \left\|\nabla^{N} \partial_{t}^{M} \varphi_{\left(k_{1}\right)}\right\|_{C_{t} L_{x}^{p}}+\left\|\nabla^{N} \partial_{t}^{M} \varphi_{\left(k_{1}\right)}^{\prime}\right\|_{C_{t} L_{x}^{p}}\lesssim r_{\|}^{\frac{1}{p}-\frac{1}{2}}\left(\frac{r_{\perp} \lambda}{r_{\|}}\right)^{N}\left(\frac{\lambda r_{\perp}\mu}{r_{\|}}\right)^{M},  \label{4.10}\tag{3.10}\\
& \left\|\nabla^{N} \psi_{(k)}\right\|_{L_{x}^{p}}+\left\|\nabla^{N} \psi_{(k)}^{\prime}\right\|_{L_{x}^{p}} \lesssim r_{\perp}^{\frac{1}{p}-\frac{1}{2}} \lambda^{N}. \label{4.11}\tag{3.11}
\end{align*}
Furthermore, it holds that
\begin{align*}
& \left\|\nabla^{N} \partial_{t}^{M} \mathbf{W}_{k}\right\|_{C_{t} L_{x}^{p}}+\frac{r_{\|}}{r_{\perp}}\left\|\nabla^{N} \partial_{t}^{M} \mathbf{W}_{k}^{c}\right\|_{C_{t} L_{x}^{p}}+r_{\perp}^{-1}\left\|\nabla^{N} \partial_{t}^{M} \boldsymbol{\Psi}_{k}\right\|_{C_{t} L_{x}^{p}} \\
& \quad \lesssim (r_{\perp}r_{\|})^{\frac{1}{p}-\frac{1}{2}}\lambda^{N}\left(\frac{\lambda r_{\perp}\mu}{r_{\|}}\right)^{M}, \label{4.12} \tag{3.12}
\end{align*}
where the implicit constants $C$ are independent of $r_{\perp}, r_{\|}, \lambda, \mu$, and $k \in \Lambda$.
\end{Lemma}

\subsection{Temporal Building Blocks.}\ \ \ \
Next, we construct the temporal building blocks that possess two crucial parameters, $\tau$ and $\sigma$, which govern concentration and oscillation over time, respectively. It is crucial to consider suitable temporal shifts $t_{k}$ within these building blocks to ensure that the supports of different temporal building blocks are disjoint. This construction is motivated by recent advancements in temporal intermittency, as detailed in references \cite{AC21,AC22,ACL22,ACX22}, which permits to obtain sharp non-uniqueness results.

More precisely, we define $\left\{g_k\right\}_{k \in \Lambda} \subset C_c^{\infty}([0, T])$ as cut-off functions with mean zero,  such that for $k \neq k^{\prime}$, the temporal supports of $g_k$ and $g_{k^{\prime}}$ are disjoint.  Furthermore, for all $k \in \Lambda$, we have
\begin{align*}
\fint_0^T g_k^2(t) \mathrm{d} t=1.
\label{4.13}\tag{3.13}
\end{align*}
Since there are finitely many wavevectors in $\Lambda$, the existence of such $\left\{g_k\right\}_{k \in \Lambda}$ can be guaranteed by $g_k=g\left(t-t_k\right)$, where $g \in C_c^{\infty}([0, T])$ has a very small support, and $\left\{t_k\right\}_{k \in \Lambda}$ are temporal shifts chosen such that the supports of $\left\{g_k\right\}_{k \in \Lambda}$ are disjoint. Subsequently, we rescale $g_k$ by
$g_{k, \tau}(t)=\tau^{\frac{1}{2}} g_k(\tau t),$ where the concentration parameter $\tau$ is defined by \eqref{4.1}. By an abuse of notation, we consider $g_{k, \tau}$ as a periodic function on $[0, T]$.
Additionally, we define
\begin{align*}
h_{k, \tau}(t):=\int_0^t\left(g_{k, \tau}^2(s)-1\right) d s, \quad t \in[0, T],
\label{4.14}\tag{3.14}
\end{align*}
and set
\begin{align*}
g_{(k)}(t):=g_{k, \tau}(\sigma t), \quad h_{(k)}(t):=h_{k, \tau}(\sigma t) .
\label{4.15}\tag{3.15}
\end{align*}

Subsequently, the function $h_{(k)}$ satisfies
\begin{align*}
\partial_t\left(\sigma^{-1} h_{(k)}\right)=g_{(k)}^2-1=g_{(k)}^2-\fint_0^T g_{(k)}^2(t) \mathrm{d} t,
\label{4.16}\tag{3.16}
\end{align*}
where $\sigma$ is given by \eqref{4.1}.

We draw upon the crucial estimates of $g_{(k)}$ and $h_{(k)}$ established in previous works \cite{AC21,YLZ22} and present them in the following lemma.

\begin{Lemma}\label{D.3} (Estimates of temporal intermittency). For $\gamma \in[1, \infty], M \in \mathbb{N}$, we have
\begin{equation*}
\left\|\partial_{t}^{M} g_{(k)}\right\|_{L_{t}^{\gamma}} \lesssim \sigma^{M} \tau^{M+\frac{1}{2}-\frac{1}{\gamma}}, \label{4.17}\tag{3.17}
\end{equation*}
where the implicit constants do not rely on both $\sigma$ and $\tau$. Moreover, we have the following estimate for $h_{(k)}$:
\begin{equation*}
\left\|h_{(k)}\right\|_{C_{t}} \leq 1. \label{4.18}\tag{3.18}
\end{equation*}
\end{Lemma}

\begin{Remark}\label{D.2} After defining the temporal concentrated function $g_{(k)}$ and 2D intermittent jets $\mathbf{W}_{k}$, we will utilize $g_{(k)} \mathbf{W}_{k}$ as the ``Building block" for Convex integration scheme. Consequently, the supports of all $g_{(k)} \mathbf{W}_{k}$ are mutually disjoint if $k\neq k'$ on $\mathbb{T}^2 \times[0,T]$.
\end{Remark}

\subsection{Velocity Perturbations.}\label{C.3.3}\ \ \ \
In the following sections, we are committed to constructing the velocity perturbations, containing the principal perturbation, incompressible corrector, and two temporal correctors. Firstly, we will provide the specific form of the amplitudes of these perturbations, which has been presented in previous research works \cite{TBV19,YL22,LTP20}.

\noindent {\bf  Amplitudes Of Perturbations.}
Choosing $\mathscr{X} :[0, \infty) \rightarrow \mathbb{R}$ to be a smooth cut-off function that satisfies
\[
\mathscr{X} (z)= \begin{cases}1, & 0 \leq z \leq 1,  \label{4.19}\tag{3.19}\\ z, & z \geq 2,\end{cases}
\]
and
\begin{equation*}
\frac{1}{2} z \leq \mathscr{X} (z) \leq 2 z, \quad \text { for } \quad z \in(1,2). \label{4.20}\tag{3.20}
\end{equation*}
Define
\begin{equation*}
\rho_{u}(t, x):=2 C_{\mathsf{R}}^{-1} \delta_{q+1} \mathscr{X} \left(\frac{\left|\mathring{\mathsf{R}}_{\ell}^{*}(t, x)\right|}{ \delta_{q+1}}\right), \label{4.21}\tag{3.21}
\end{equation*}
where $C_{\mathsf{R}}$ is the positive constant defined in the Geometric Lemma \ref{G.1}. By \eqref{4.19}, \eqref{4.20} and \eqref{4.21},
\begin{equation*}
\left|\frac{\mathring{\mathsf{R}}_{\ell}^{*}}{\rho_{u}}\right|=\left|\frac{\mathring{\mathsf{R}}_{\ell}^{*}}{2 C_{\mathsf{R}}^{-1}  \delta_{q+1} \mathscr{X} \left( \delta_{q+1}^{-1}\left|\mathring{\mathsf{R}}_{\ell}^{*}\right|\right)}\right| \leq C_{\mathsf{R}}, \label{4.22}\tag{3.22}
\end{equation*}
and for any $p \in[1, \infty]$,
\begin{align*}
& \rho_{u} \geq C_{\mathsf{R}}^{-1} \delta_{q+1},\label{4.23}\tag{3.23}
\end{align*}
\begin{align*}
\|\rho_{u}\|_{L_{t, x}^p} \lesssim C_{\mathsf{R}}^{-1}\left( \delta_{q+1}+\left\|\mathring{\mathsf{R}}_{\ell}^{*}\right\|_{L_{t, x}^p}\right).\label{4.24}\tag{3.24}
\end{align*}
Furthermore, combining \eqref{222}, \eqref{4.23} with the standard H\"{o}lder estimates (see \cite{TB15}), for $1 \leq N \leq 5$,
\begin{gather*}
\|\rho_{u}\|_{C_{t, x}} \lesssim \ell^{-2}, \quad\|\rho_{u}\|_{C_{t, x}^{N}} \lesssim \ell^{-2 N} \label{4.25}\tag{3.25}\\
\left\|\rho_{u}^{1 / 2}\right\|_{C_{t, x}} \lesssim \ell^{-1}, \quad\left\|\rho_{u}^{1 / 2}\right\|_{C_{t, x}^{N}} \lesssim \ell^{-2 N},  \label{4.26}\tag{3.26}\\
\left\|\rho_{u}^{-1}\right\|_{C_{t, x}} \lesssim \ell^{-1}, \quad\left\|\rho_{u}^{-1}\right\|_{C_{t, x}^{N}} \lesssim \ell^{-2 N}, \label{4.27}\tag{3.27}
\end{gather*}
where the implicit constants are independent of the variable $q$. To ensure temporal compatibility between the perturbations and the Reynolds stress $\mathring{\mathsf{R}}_{\ell}^{*}$ discussed in Section \ref{B}, we employ a smoothly temporal cut-off function $f_{u}:[0, T] \rightarrow[0,1]$ that fulfills
\begin{itemize}
  \item $0 \leq f_{u} \leq 1$ and $f \equiv 1$ on $\operatorname{supp}_{t} \mathring{\mathsf{R}}_{\ell}^{*} ;$
  \item $\operatorname{supp}_t f_{u} \subseteq N_{\ell}\left(\operatorname{supp}_t \mathring{\mathsf{R}}_{\ell}^{*}\right);$
  \item $\|f_{u}\|_{C_{t}^{N}} \lesssim \ell^{-N}, \quad 1 \leq N \leq 5$.
\end{itemize}

Now, we provide the specific form of amplitude as follows
\begin{equation*}
a_{(k)}(t, x):=\rho_{u}^{\frac{1}{2}}(t, x) f_{u}(t) \gamma_{(k)}\left(\operatorname{Id}-\frac{\mathring{\mathsf{R}}_{\ell}^{*}(t, x)}{\rho_{u}(t, x)}\right), \quad k \in \Lambda, \label{4.28}\tag{3.28}
\end{equation*}
where the symbols $\gamma_{(k)}$ and $\Lambda$ are given in the Geometric Lemma \ref{G.1}. Furthermore, we give the analytic estimates for the amplitudes, the derivation of which follows a similar proof in \cite{YLZ22}.
\begin{Lemma}\label{D.5}(Estimates of amplitudes)
For $1 \leq N \leq 5, k \in \Lambda$, it holds
\begin{align*}
\left\|a_{(k)}\right\|_{L_{t, x}^{2}} & \lesssim \delta_{q+1}^{\frac{1}{2}},  \label{4.29}\tag{3.29}\\
\left\|a_{(k)}\right\|_{C_{t, x}} & \lesssim \ell^{-1}, \quad\left\|a_{(k)}\right\|_{C_{t, x}^{N}} \lesssim \ell^{-4 N} , \label{4.30}\tag{3.30}
\end{align*}
where the implicit constants are independent of the variable $q$.
\end{Lemma}
\noindent {\bf Velocity Perturbations.}
After the selection of ``Building block'' and the introduction of amplitudes, we now begin to construct the velocity perturbations. Firstly, the principal part of the velocity perturbations $w_{q+1}^{(p)}$, is defined as
\begin{equation*}
w_{q+1}^{(p)}:=\sum_{k \in \Lambda} a_{(k)} g_{(k)} \mathbf{W}_{k}. \label{4.31}\tag{3.31}
\end{equation*}

By applying Geometric Lemma \ref{G.1} and \eqref{4.28}, the concentrated Reynolds stress $\mathring{\mathsf{R}}_{\ell}^{*}$ can be eliminated by the zero frequency part of $w_{q+1}^{(p)} \otimes w_{q+1}^{(p)}$, i.e.,
\begin{align*}
w_{q+1}^{(p)} \otimes w_{q+1}^{(p)}+\mathring{\mathsf{R}}_{\ell}^{*}&=\rho_{u} f_{u}^2 \operatorname{Id}+\sum_{k \in \Lambda} a_{(k)}^2 g_{(k)}^2 \mathbb{P}_{\neq 0}\left(\mathbf{W}_{k} \otimes \mathbf{W}_{k}\right)\\
& +\sum_{k \in \Lambda} a_{(k)}^{2}\left(g_{(k)}^{2}-1\right) \fint_{\mathbb{T}^{2}} \mathbf{W}_{k} \otimes \mathbf{W}_{k} \mathrm{d} x .\label{4.32}\tag{3.32}
\end{align*}
In the given context,  $\mathbb{P}_{\neq 0}$ represents the spatial projection onto the set of nonzero Fourier modes. Since $w_{q+1}^{(p)}$ lacks divergence-free property, we introduce the incompressibility corrector defined as:
\begin{equation*}
w_{q+1}^{(c)}:=\sum_{k \in \Lambda} a_{(k)} g_{(k)}{\bf W}_{k}^{c}+(\lambda r_{\perp}N_{\Lambda})^{-1}\nabla^{\bot} a_{(k)} g_{(k)}{\bf\Psi}_{k}, \label{4.33}\tag{3.33}
\end{equation*}
where $\mathbf{W}_{k}^{c}$ and ${\bf\Psi}_{k} $ are respectively defined by \eqref{4.6} and \eqref{4.7}. subsequently,
\begin{align*}
w_{q+1}^{(p)}+w_{q+1}^{(c)}&=\sum_{k \in \Lambda}a_{(k)} g_{(k)}({\bf W}_{k}+{\bf W}_{k}^{c})
+(\lambda r_{\perp}N_{\Lambda})^{-1}\nabla^{\bot} a_{(k)} g_{(k)}{\bf\Psi}_{k}\\
&=(\lambda r_{\perp}N_{\Lambda})^{-1}\sum_{k \in \Lambda} g_{(k)}[a_{(k)}\nabla^{\bot}{\bf\Psi}_{k}
+\nabla^{\bot} a_{(k)}{\bf\Psi}_{k}]\\
&= (\lambda r_{\perp}N_{\Lambda})^{-1}\sum_{k \in \Lambda}\nabla^{\bot}[a_{(k)} g_{(k)}{\bf\Psi}_{k}],\label{4.34}\tag{3.34}
\end{align*}
and thus
\begin{equation*}
\operatorname{div}\left(w_{q+1}^{(p)}+w_{q+1}^{(c)}\right)=0. \label{4.35}\tag{3.35}
\end{equation*}

To address the high-frequency spatial and temporal errors in \eqref{4.32}, we introduce two additional  temporal correctors. The first temporal corrector $w_{q+1}^{(t)}$ is defined as
\begin{equation*}
w_{q+1}^{(t)}:=-\mu^{-1} \sum_{k \in \Lambda} \mathbb{P}_{H} \mathbb{P}_{\neq 0}\left(a_{(k)}^{2} g_{(k)}^{2} \varphi_{\left(k_{1}\right)} ^{2} (\psi_{(k)}^{\prime})^{2} k_{1}\right) \label{4.36}\tag{3.36}
\end{equation*}
to balance the high spatial frequency oscillations:
\begin{align*}
\partial_{t} w_{q+1}^{(t)}&+\sum_{k \in \Lambda} \mathbb{P}_{\neq 0}\left(a_{(k)}^{2} g_{(k)}^{2} \operatorname{div}\left(\mathbf{W}_{k} \otimes \mathbf{W}_{k}\right)\right)=\\
 & \left(\nabla \Delta^{-1} \operatorname{div}\right) \mu^{-1} \sum_{k \in \Lambda} \mathbb{P}_{\neq 0} \partial_{t}\left(a_{(k)}^{2} g_{(k)}^{2} \varphi_{\left(k_{1}\right)} ^{2} (\psi_{(k)}^{\prime})^{2} k_{1}\right) \\
& -\mu^{-1} \sum_{k \in \Lambda} \mathbb{P}_{\neq 0}\left(\partial_{t}\left(a_{(k)}^{2} g_{(k)}^{2}\right) \varphi_{\left(k_{1}\right)} ^{2} (\psi_{(k)}^{\prime})^{2} k_{1}\right) , \label{4.37}\tag{3.37}
\end{align*}
where the first term on the right-hand-side can be processed by Helmholtz-Leray projector $\mathbb{P}_{H}$.

The other temporal corrector $w_{q+1}^{(o)}$ is defined by
\begin{equation*}
w_{q+1}^{(o)}:=-\sigma^{-1} \sum_{k \in \Lambda} \mathbb{P}_{H} \mathbb{P}_{\neq 0}\left(h_{(k)} \fint_{\mathbb{T}^{2}} \mathbf{W}_{k} \otimes \mathbf{W}_{k} \mathrm{d} x \nabla\left(a_{(k)}^{2}\right)\right) \label{4.38}\tag{3.38}
\end{equation*}
to balance the high temporal frequency oscillations:
\begin{align*}
& \partial_{t} w_{q+1}^{(o)}+\sum_{k \in \Lambda} \mathbb{P}_{\neq 0}\left(\left(g_{(k)}^{2}-1\right) \fint_{\mathbb{T}^{2}} \mathbf{W}_{k} \otimes \mathbf{W}_{k} \mathrm{d} x \nabla\left(a_{(k)}^{2}\right)\right) \\
= & \left(\nabla \Delta^{-1} \operatorname{div}\right) \sigma^{-1} \sum_{k \in \Lambda} \mathbb{P}_{\neq 0} \partial_{t}\left(h_{(k)} \fint_{\mathbb{T}^{2}} \mathbf{W}_{k} \otimes \mathbf{W}_{k} \mathrm{d} x \nabla\left(a_{(k)}^{2}\right)\right) \\
& -\sigma^{-1} \sum_{k \in \Lambda} \mathbb{P}_{\neq 0}\left(h_{(k)} \fint_{\mathbb{T}^{2}} \mathbf{W}_{k} \otimes \mathbf{W}_{k} \mathrm{d} x \partial_{t} \nabla\left(a_{(k)}^{2}\right)\right) . \label{4.39}\tag{3.39}
\end{align*}
where the right-hand-side above only remain the low frequency part $\partial_{t} \nabla\left(a_{(k)}^{2}\right)$ and the harmless pressure term.

In summary, the velocity perturbation $w_{q+1}$ at level $q+1$ is defined as the sum of four components:
\begin{equation*}
w_{q+1}:=w_{q+1}^{(p)}+w_{q+1}^{(c)}+w_{q+1}^{(t)}+w_{q+1}^{(o)} . \label{4.40}\tag{3.40}
\end{equation*}
It is easy to verify that $w_{q+1}$ are both mean-free and divergence-free. Subsequently, the velocity field at level $q+1$ is determined by
\begin{equation*}
u_{q+1}:=u_{\ell}+w_{q+1}, \label{4.41}\tag{3.41}
\end{equation*}
where $u_{\ell}$ is defined as in \eqref{212}.

Furthermore, we present the crucial estimates associated with these velocity perturbations in the following Lemma \ref{D.6}.
\begin{Lemma}\label{D.6} (Estimates of velocity perturbations) Given any $p\in(1, \infty), \gamma \in[1, \infty]$ and integers $0 \leq N \leq 5$, it holds that
\begin{align*}
&\left\|\nabla^{N} w_{q+1}^{(p)}\right\|_{L_{t}^{\gamma} L_{x}^{p}} \lesssim \ell^{-1} \lambda^{N} (r_{\perp} r_{\|})^{\frac{1}{p}-\frac{1}{2}} \tau^{\frac{1}{2}-\frac{1}{\gamma}},  \label{4.42}\tag{3.42}\\
&\left\|\nabla^{N} w_{q+1}^{(c)}\right\|_{L_{t}^{\gamma} L_{x}^{p}} \lesssim \ell^{-6} \lambda^{N} (r_{\perp} r_{\|})^{\frac{1}{p}-\frac{1}{2}} \tau^{\frac{1}{2}-\frac{1}{\gamma}}  \frac{r_{\perp}}{ r_{\|}},\label{4.43}\tag{3.43}\\
&\left\|\nabla^{N} w_{q+1}^{(t)}\right\|_{L_{t}^{\gamma} L_{x}^{p}} \lesssim \ell^{-2} \lambda^{N} \mu^{-1} (r_{\perp}r_{\|})^{\frac{1}{p}-1} \tau^{1-\frac{1}{\gamma}},  \label{4.44}\tag{3.44}\\
&\left\|\nabla^{N} w_{q+1}^{(o)}\right\|_{L_{t}^{\gamma} L_{x}^{p}} \lesssim \ell^{-4N-5} \sigma^{-1}. \label{4.45}\tag{3.45}
\end{align*}
Especially, for integers $1 \leq N \leq 5$, we have
\begin{align*}
\left\|w_{q+1}^{(p)}\right\|_{C_{t, x}^N}+\left\|w_{q+1}^{(c)}\right\|_{C_{t, x}^N}+\left\|w_{q+1}^{(t)}\right\|_{C_{t, x}^N}+\left\|w_{q+1}^{(o)}\right\|_{C_{t, x}^N} \lesssim \lambda^{3N+3}, \label{4.47}\tag{3.46}
\end{align*}
where the implicit constants are independent of $\lambda$.
\end{Lemma}
{\bf Proof.} By  \eqref{4.12}, \eqref{4.17}, \eqref{4.32}, \eqref{4.33} and Lemma \ref{D.5}, for any $p \in(1, \infty)$,
$$
\begin{aligned}
\left\|\nabla^{N} w_{q+1}^{(p)}\right\|_{L_{t}^{\gamma} L_{x}^{p}} & \lesssim \sum_{k \in \Lambda} \sum_{N_{1}+N_{2}=N}\left\|a_{(k)}\right\|_{C_{t, x}^{N_{1}}}\left\|g_{(k)}\right\|_{L_{t}^{\gamma}}\left\|\nabla^{N_{2}} \mathbf{W}_{k}\right\|_{C_{t} L_{x}^{p}} \\
& \lesssim \ell^{-1} \lambda^{N} (r_{\perp} r_{\|})^{\frac{1}{p}-\frac{1}{2}} \tau^{\frac{1}{2}-\frac{1}{\gamma}},
\end{aligned}
$$
and
$$
\begin{aligned}
&\left\|\nabla^{N} w_{q+1}^{(c)}\right\|_{L_{t}^{\gamma} L_{x}^{p}} \\
& \lesssim \sum_{k \in \Lambda}\left\|g_{(k)}\right\|_{L_{t}^{\gamma}} \sum_{N_{1}+N_{2}=N}\left(\left\|a_{(k)}\right\|_{C_{t, x}^{N_{1}}}\left\|\nabla^{N_{2}} \mathbf{W}_{k}^{c}\right\|_{C_{t} L_{x}^{p}}+(\lambda r_{\perp})^{-1}\left\|a_{(k)}\right\|_{C_{t, x}^{N_{1}+1}}\left\|\nabla^{N_{2}} \boldsymbol{\Psi}_{k} \right\|_{C_{t} L_{x}^{p}}\right) \\
& \lesssim \sum_{N_{1}+N_{2}=N} \tau^{\frac{1}{2}-\frac{1}{\gamma}}\left(\ell^{-4 N_{1}-1} \lambda^{N_{2}} \frac{r_{\perp}}{r_{\|}}(r_{\perp} r_{\|})^{\frac{1}{p}-\frac{1}{2}}+(\lambda r_{\perp})^{-1}\ell^{-4 N_{1}-6} \lambda^{N_{2}} r_{\perp} ( r_{\perp}r_{\|})^{\frac{1}{p}-\frac{1}{2}}\right) \\
& \lesssim \tau^{\frac{1}{2}-\frac{1}{\gamma}}\left(\ell^{-1} \lambda^{N} \frac{r_{\perp}}{r_{\|}}(r_{\perp} r_{\|})^{\frac{1}{p}-\frac{1}{2}}+\ell^{-6} \lambda^{N-1}(r_{\perp}r_{\|})^{\frac{1}{p}-\frac{1}{2}}\right) \\
& \lesssim \ell^{-6} \lambda^{N} (r_{\perp}r_{\|})^{\frac{1}{p}-\frac{1}{2}} \tau^{\frac{1}{2}-\frac{1}{\gamma}}\frac{r_{\perp}}{r_{\|}} ,
\end{aligned}
$$
then, the \eqref{4.42}, \eqref{4.43} are verified.

Regarding the temporal correctors, based on \eqref{4.36} and \eqref{4.38}, Lemmas \ref{D.1}, \ref{D.3}, and \ref{D.5}, as well as the boundedness of the operators $\mathbb{P}{\neq 0}$ and $\mathbb{P}_{H}$ in the space $L_{x}^{p}$, we deduce that
$$
\begin{aligned}
&\left\|\nabla^{N} w_{q+1}^{(t)}\right\|_{L_{t}^{\gamma} L_{x}^{p}}\\
& \lesssim \mu^{-1} \sum_{k \in \Lambda}\left\|g_{(k)}^{2}\right\|_{L_{t}^{\gamma}} \sum_{N_{1}+N_{2}+N_{3}=N}\left\|\nabla^{N_{1}}\left(a_{(k)}^{2}\right)\right\|_{C_{t, x}}\left\|\nabla^{N_{2}}\left(\varphi_{\left(k_{1}\right)}^{2}\right)\right\|_{C_{t} L_{x}^{p}}\left\|\nabla^{N_{3}}\left(\psi_{(k)}^{\prime}\right)^{2}\right\|_{L_{x}^{p}} \\
& \lesssim \mu^{-1} \tau^{1-\frac{1}{\gamma}} \sum_{N_{1}+N_{2}+N_{3}=N} \ell^{-4 N_{1}-2} (\lambda\frac{r_{\perp}}{r_{\|}})^{N_{2}} \lambda^{N_{3}} (r_{\|}r_{\perp})^{\frac{1}{p}-1} \\
& \lesssim \ell^{-2}\mu^{-1} \lambda^{N} (r_{\|}r_{\perp})^{\frac{1}{p}-1} \tau^{1-\frac{1}{\gamma}},
\end{aligned}
$$
and
$$
\left\|\nabla^{N} w_{q+1}^{(o)}\right\|_{L_{t}^{\gamma} L_{x}^{p}} \lesssim \sigma^{-1} \sum_{k \in \Lambda}\left\|h_{(k)}\right\|_{C_{t}}\left\|\nabla^{N+1}\left(a_{(k)}^{2}\right)\right\|_{C_{t, x}} \lesssim \ell^{-4N-5} \sigma^{-1}
$$
which yields \eqref{4.44} and \eqref{4.45}.

It remains to prove the $C_{t,x}^N$-estimate \eqref{4.47} of perturbations. By Lemmas \ref{D.1}, \ref{D.3}, and \ref{D.5},

\begin{align*}
\left\|w_{q+1}^{(p)}\right\|_{C_{t, x}^N} & \lesssim \sum_{k \in \Lambda}\left\|a_{(k)}\right\|_{C_{t, x}^N} \sum_{0 \leq N_1+N_2 \leq N}\left\|g_{(k)}\right\|_{C_t^{N_1}}\left\|\mathbf{W}_{(k)}\right\|_{C_{t, x}^{N_2}} \\
& \lesssim \sum_{0 \leq N_1+N_2 \leq N} \sum_{N_{21}+N_{22}=N_2} \ell^{-4 N-1} \sigma^{N_1} \tau^{N_1+\frac{1}{2}} r_{\perp}^{-\frac{1}{2}} r_{\|}^{-\frac{1}{2}} \lambda^{N_{21}}\left(\frac{ \lambda r_{\perp} \mu}{r_{\|}}\right)^{N_{22}} \\
& \lesssim \lambda^{3 N+3}, \label{347}\tag{3.47}
\end{align*}
where we also used \eqref{2.3} and \eqref{4.1} in the last step.

Similarly, we have
\begin{align*}
\left\|w_{q+1}^{(c)}\right\|_{C_{t, x}^N}&\lesssim \sum_{k \in \Lambda}\left\|a_{(k)}\right\|_{C_{t, x}^{N+1}} \sum_{0 \leq N_1+N_2 \leq N}\left\|g_{(k)}\right\|_{C_t^{N_1}}\left(\left\|W_{(k)}^c\right\|_{C_{t, x}^{N_2}}+(\lambda r_{\perp})^{-1}\left\|\mathbf{\Psi}_{(k)}\right\|_{C_{t, x}^{N_2}}\right) \\
& \lesssim \sum_{0 \leq N_1+N_2 \leq N} \ell^{-4 N-4} \sigma^{N_1} \tau^{N_1+\frac{1}{2}} r_{\perp}^{-\frac{1}{2}} r_{\|}^{-\frac{1}{2}} \lambda^{N_{21}}\left(\frac{ \lambda r_{\perp} \mu}{r_{\|}}\right)^{N_{2}}(\frac{r_{\perp} }{r_{\|}}+\lambda^{-1}) \\
& \lesssim \lambda^{3 N+3}.
\label{4.50}\tag{3.48}
\end{align*}

Moreover, by Sobolev's embedding $W^{1,4}\left(\mathbb{T}^2\right) \hookrightarrow L^{\infty}\left(\mathbb{T}^2\right)$ and the boudedness of operators $\mathbb{P}_H \mathbb{P}_{\neq 0}$ in the space $L_x^4$,

\begin{align*}
\left\|w_{q+1}^{(t)}\right\|_{C_{t, x}^N} & \lesssim \mu^{-1} \sum_{k \in \Lambda}\left\|a_{(k)}^2 g_{(k)}^2 \varphi_{\left(k_1\right)}^2 (\psi_{(k)}^{\prime})^2\right\|_{C_t^N W_x^{N+1,4}} \\
& \lesssim \mu^{-1} \sum_{k \in \Lambda} \sum_{0 \leq N_1+N_2 \leq N+1}\left\|a_{(k)}^2\right\|_{C_{t, x}^N}\left\|g_{(k)}^2 \varphi_{\left(k_1\right)}^2 (\psi_{(k)}^{\prime})^2\right\|_{C_t^{N_1} C_x^{N_2}} \\
& \leq\lambda^{3 N+2},\label{349}\tag{3.49}
\end{align*}
where we also used Lemmas \ref{D.1}, \ref{D.3}, and \ref{D.5} in the last step.
Finally, arguing as above we get
\begin{align*}
\left\|w_{q+1}^{(o)}\right\|_{C_{t, x}^N} & \lesssim \sigma^{-1} \sum_{k \in \Lambda_u \cup \Lambda_B}\left\|h_{(k)} \nabla\left(a_{(k)}^2\right)\right\|_{C_t^N W_x^{N+1,4}} \\
& \lesssim \sigma^{-1} \sum_{k \in \Lambda}\left\|h_{(k)}\right\|_{C_t^{N}}\left\|\nabla\left(a_{(k)}^2\right)\right\|_{C_{t, x}^{N+1}} \\
& \lesssim \sigma^{N-1} \tau^{N} \ell^{-4(N+2)-1} \\
&\lesssim \lambda^{3N+1},\label{350}\tag{3.50}
\end{align*}
where the last two steps were due to \eqref{2.3} and \eqref{4.1}.
Combining \eqref{347}-\eqref{350} altogether and using \eqref{4.1} we conclude that
\begin{align*}
\left\|w_{q+1}^{(p)}\right\|_{C_{t, x}^N}+\left\|w_{q+1}^{(c)}\right\|_{C_{t, x}^N}+\left\|w_{q+1}^{(t)}\right\|_{C_{t, x}^N}+\left\|w_{q+1}^{(o)}\right\|_{C_{t, x}^N} \leq \lambda^{3N+3}. \label{351}\tag{3.51}
\end{align*}

Consequently, the proof of Lemma \ref{D.6} is finished.

\subsection{Inductive Estimates For Velocity Perturbations.}\ \ \ \
In this section, we aim to verify the iterative estimates \eqref{2.6}, and \eqref{2.11}-\eqref{2.13} pertaining to the velocity perturbations. We invoke the $L^{p}$ de-correlation Lemma \ref{H.2}, substituting $f$ with $a_{(k)}$ and $g$ with $g_{(k)} \varphi_{\left(k_{1}\right)} \psi_{(k)}^{\prime}$ while setting $\sigma=\lambda^{2 \varepsilon}$. Subsequently, we utilize \eqref{2.2}-\eqref{2.4} alongside Lemmas \ref{D.1}, \ref{D.3}, and \ref{D.5} to obtain:
\begin{gather*}
\left\|w_{q+1}^{(p)}\right\|_{L_{t, x}^{2}} \lesssim \sum_{k \in \Lambda}\left(\left\|a_{(k)}\right\|_{L_{t, x}^{2}}\left\|g_{(k)}\right\|_{L_{t}^{2}}\left\|\varphi_{\left(k_{1}\right)} \psi_{(k)}^{\prime}\right\|_{C_{t} L_{x}^{2}}\right. \\
\left.+\sigma^{-\frac{1}{2}}\left\|a_{(k)}\right\|_{C_{t, x}^{1}}\left\|g_{(k)}\right\|_{L_{t}^{2}}\left\|\varphi_{\left(k_{1}\right)} \psi_{(k)}^{\prime}\right\|_{C_{t} L_{x}^{2}}\right) \\
\lesssim \delta_{q+1}^{\frac{1}{2}}+\ell^{-4} \lambda_{q+1}^{-\varepsilon} \lesssim \delta_{q+1}^{\frac{1}{2}} . \label{4.53}\tag{3.52}
\end{gather*}

In light of \eqref{2.3}, by employing \eqref{4.53} and invoking Lemma \ref{D.6}, we establish an upper bound for the velocity perturbation,
\begin{align*}
\left\|w_{q+1}\right\|_{L_{t, x}^{2}} & \lesssim\left\|w_{q+1}^{(p)}\right\|_{L_{t, x}^{2}}+\left\|w_{q+1}^{(c)}\right\|_{L_{t, x}^{2}}+\left\|w_{q+1}^{(t)}\right\|_{L_{t, x}^{2}}+\left\|w_{q+1}^{(o)}\right\|_{L_{t, x}^{2}} \\
& \lesssim \delta_{q+1}^{\frac{1}{2}}+\ell^{-6} r_{\perp} r_{\|}^{-1}+\ell^{-2} \mu^{-1} (r_{\perp}r_{\|})^{-\frac{1}{2}}\tau^{\frac{1}{2}}+\ell^{-5} \sigma^{-1} \lesssim \delta_{q+1}^{\frac{1}{2}}, \label{4.54}\tag{3.53}
\end{align*}
and
\begin{align*}
\left\|w_{q+1}\right\|_{L_{t}^{1} L_{x}^{2}} & \lesssim\left\|w_{q+1}^{(p)}\right\|_{L_{t}^{1} L_{x}^{2}}+\left\|w_{q+1}^{(c)}\right\|_{L_{t}^{1} L_{x}^{2}}+\left\|w_{q+1}^{(t)}\right\|_{L_{t}^{1} L_{x}^{2}}+\left\|w_{q+1}^{(o)}\right\|_{L_{t}^{1} L_{x}^{2}} \label{4.55}\tag{3.54}\\
& \lesssim \ell^{-1} \tau^{-\frac{1}{2}}+\ell^{-6} \frac{r_{\perp}}{r_{\|}} \tau^{-\frac{1}{2}}+\ell^{-2} \mu^{-1} (r_{\perp}r_{\|})^{-\frac{1}{2}} +\ell^{-5} \sigma^{-1} \lesssim \lambda_{q+1}^{-\varepsilon}.
\end{align*}

Next, we proceed to verify the iterative estimates for $u_{q+1}$, we apply \eqref{2.6}, \eqref{4.41}, \eqref{4.47} to derive the following results:
$$\|u_{q+1}\|_{C_{t,x}^{1}}\leq\|u_{\ell}\|_{C_{t,x}^{1}}+
\|w_{q+1}\|_{C_{t,x}^{1}}\lesssim\lambda_{q}^{7}+\lambda_{q+1}^{6}
\lesssim \lambda_{q+1}^{7}.$$

Furthermore, given the fact $\lambda_{q}^{-13} \ll \delta_{q+2}^{1 / 2}$, we employ \eqref{2.2}, \eqref{2.3}, \eqref{220} and \eqref{4.54} to deduce that
\begin{align*}
\left\|u_{q}-u_{q+1}\right\|_{L_{t, x}^{2}} & \leq\left\|u_{q}-u_{\ell}\right\|_{L_{t, x}^{2}}+\left\|u_{\ell}-u_{q+1}\right\|_{L_{t, x}^{2}} \\
& \lesssim\left\|u_{q}-u_{\ell}\right\|_{L_{t,x}^{\infty}}+\left\|w_{q+1}\right\|_{L_{t, x}^{2}} \\
& \lesssim \lambda_{q}^{-13}+\delta_{q+1}^{\frac{1}{2}} \leq\delta_{q+1}^{\frac{1}{2}}, \label{4.58}\tag{3.55}
\end{align*}
and
\begin{align*}
\left\|u_{q}-u_{q+1}\right\|_{L_{t}^{1} L_{x}^{2}} & \lesssim\left\|u_{q}-u_{\ell}\right\|_{L_{t}^{\infty} L_{x}^{2}}+\left\|w_{q+1}\right\|_{L_{t}^{1} L_{x}^{2}} \\
& \lesssim \lambda_{q}^{-13}+\lambda_{q+1}^{-\varepsilon} \leq \delta_{q+2}^{\frac{1}{2}}, \label{4.59}\tag{3.56}
\end{align*}
where we set $a$ to be sufficiently large, ensuring the validity of the final inequalities in \eqref{4.59}.

Then, for any $\alpha \in[1 , \frac{3}{2})$, by \eqref{2.3}, \eqref{4.1} and Lemma \ref{D.6},
\begin{align*}
\left\|u_{q+1}-u_{q}\right\|_{L_{t}^{\gamma} L_{x}^{p}} & \lesssim\left\|u_{\ell}-u_{q}\right\|_{L_{t}^{\gamma} L_{x}^{p}}+\left\|w_{q+1}\right\|_{L_{t}^{\gamma} L_{x}^{p}} \\
& \lesssim \lambda_{q}^{-13}+\ell^{-1}  (r_{\perp} r_{\|})^{\frac{1}{p}-\frac{1}{2}} \tau^{\frac{1}{2}-\frac{1}{\gamma}}+\ell^{-17} \sigma^{-1} \\
& \lesssim \lambda_{q}^{-13}+\lambda_{q+1}^{2 \alpha-1-\frac{2}{p}-\frac{4 \alpha-4}{\gamma}+\varepsilon\left(2+\frac{12}{p}-\frac{16}{\gamma}\right)}+\lambda_{q+1}^{-\varepsilon}. \label{4.62}\tag{3.57}
\end{align*}
Incorporating the considerations from \eqref{2.4}, we have
\begin{equation*}
2 \alpha-1-\frac{2}{p}-\frac{4 \alpha-4}{\gamma}+\varepsilon\left(2+\frac{12}{p}-\frac{16}{\gamma}\right) \leq 2 \alpha-1-\frac{2}{p}-\frac{4 \alpha-4}{\gamma}+14 \varepsilon<-6 \varepsilon,\label{4.63} \tag{3.58}
\end{equation*}
which yields that
\begin{equation*}
\left\|u_{q+1}-u_{q}\right\|_{L_{t}^{\gamma} L_{x}^{p}} \leq \delta_{q+2}^{\frac{1}{2}}. \label{4.64}\tag{3.59}
\end{equation*}

Hence, the iteration estimates \eqref{2.6}, and \eqref{2.11}-\eqref{2.13} are verified.

{\centering
\section{Reynolds Stress}\label{D}}
The objective of this section is to establish the validity of the inductive estimates \eqref{2.10}-\eqref{2.8} for the new Reynolds stress $\mathring{\mathsf{R}}_{q+1}$ in the super-critical space \eqref{*} when $\alpha \in[1,\frac{3}{2})$.

The important role here is played by the inverse divergence operator $\mathcal{R}: C^{\infty}\left(\mathbb{T}^2, \mathbb{R}^2\right) \rightarrow C^{\infty}\left(\mathbb{T}^2, \mathcal{S}_0^{2 \times 2}\right)$ is defined as in \cite{CD13}, specifically
$$ (\mathcal{R} v)_{i j}=-\Delta^{-1} \partial_k \delta_{i j}v_k+\Delta^{-1} \partial_i \delta_{j k}v_k+\Delta^{-1} \partial_j \delta_{i k} v_k.
$$
By a direct computation, one can also show that
$\dive(\mathcal{R}v)=v$, where $v$ is a  mean-free smooth function satisfying $\int_{\mathbb{T}^{2}} v dx=0$. The antidivergence operator $\mathcal{R}$ is bounded on $L^p\left(\mathbb{T}^2\right)$ for any $1 < p < \infty$ (see also \cite{MS18}).\\
\subsection{Decomposition Of Reynolds Stress.}\ \ \ \
From \eqref{2.1} and \eqref{4.41}, we deduce that the new Reynolds stress $\mathring{\mathsf{R}}_{q+1}$ is given as follows:
\begin{align*}
& \operatorname{div} \mathring{\mathsf{R}}_{q+1}-\nabla(\mathsf{P}_{q+1}-\mathsf{P}_{\ell})=\underbrace{\partial_{t}\left(w_{q+1}^{(p)}+w_{q+1}^{(c)}\right)+(-\Delta)^{\alpha} w_{q+1}+\operatorname{div}\left(u_{\ell} \otimes w_{q+1}+w_{q+1} \otimes u_{\ell}\right)}_{\operatorname{div} \mathring{\mathsf{R}}_{\{l i n\}}+\nabla\mathsf{P}_{\{l i n\}}} \\
& +\underbrace{\operatorname{div}\left(w_{q+1}^{(p)} \otimes w_{q+1}^{(p)}+\mathring{\mathsf{R}}_{\ell}^{*}\right)+\partial_{t} w_{q+1}^{(t)}+\partial_{t} w_{q+1}^{(o)}}_{\operatorname{div} \mathring{\mathsf{R}}_{\{osc\}}+\nabla\mathsf{P}_{\{o s c\}}}\label{5.1}\tag{4.1} \\
& +\underbrace{\operatorname{div}\left(\left(w_{q+1}^{(c)}+w_{q+1}^{(t)}+w_{q+1}^{(o)}\right) \otimes w_{q+1}+w_{q+1}^{(p)} \otimes\left(w_{q+1}^{(c)}+w_{q+1}^{(t)}+w_{q+1}^{(o)}\right)\right)}_{\operatorname{div} \mathring{\mathsf{R}}_{\{cor\}}+\nabla\mathsf{P}_{\{c o r\}}}\\
&+\underbrace{\operatorname{div}\left(u_{\ell} \mathring{\otimes} u_{\ell}-\left(u_q \mathring{\otimes} u_q\right) *_x \Phi_{\ell} *_t \widetilde{\Phi}_{\ell}\right)}_{\text {div } \mathring{\mathsf{R}}_{\{com\}}}.
\end{align*}

By employing the inverse divergence operator designated as $\mathcal{R}$, the definition of Reynolds stress at the level $q+1$ becomes:
\begin{equation*}
\mathring{\mathsf{R}}_{q+1}:=\mathring{\mathsf{R}}_{\{l i n\}}+\mathring{\mathsf{R}}_{\{osc\}}+\mathring{\mathsf{R}}_{\{cor\}}
+\mathring{\mathsf{R}}_{\{com\}}, \label{5.2}\tag{4.2}
\end{equation*}
where the linear error
\begin{equation*}
\mathring{\mathsf{R}}_{ \{lin \}}:=\mathcal{R}\left(\partial_{t}\left(w_{q+1}^{(p)}+w_{q+1}^{(c)}\right)\right)+\mathcal{R}(-\Delta)^{\alpha} w_{q+1}+\mathcal{R} \mathbb{P}_{H} \operatorname{div}\left(u_{\ell} \mathring{\otimes} w_{q+1}+w_{q+1} \mathring{\otimes} u_{\ell}\right) , \label{5.3}\tag{4.3}
\end{equation*}
the oscillation error
\begin{align*}
\mathring{\mathsf{R}}_{\{osc\}}:= & \sum_{k \in \Lambda} \mathcal{R} \mathbb{P}_{H} \mathbb{P}_{\neq 0}\left(g_{(k)}^{2} \mathbb{P}_{\neq 0}\left(\mathbf{W}_{k} \otimes \mathbf{W}_{k}\right) \nabla\left(a_{(k)}^{2}\right)\right) \\
& -\mu^{-1} \sum_{k \in \Lambda} \mathcal{R} \mathbb{P}_{H} \mathbb{P}_{\neq 0}\left(\partial_{t}\left(a_{(k)}^{2} g_{(k)}\right) \varphi_{\left(k_{1}\right)}^{2} (\psi_{(k)}^{\prime})^{2} k_{1}\right) \\
& -\sigma^{-1} \sum_{k \in \Lambda} \mathcal{R} \mathbb{P}_{H} \mathbb{P}_{\neq 0}\left(h_{(k)} \fint_{\mathbb{T}^{2}} \mathbf{W}_{k} \otimes \mathbf{W}_{k} \mathrm{d} x \partial_{t} \nabla\left(a_{(k)}^{2}\right)\right) , \label{5.4}\tag{4.4}
\end{align*}
the corrector error
\begin{equation*}
\mathring{\mathsf{R}}_{\{cor\}}:=\mathcal{R} \mathbb{P}_{H} \operatorname{div}\left(w_{q+1}^{(p)} \mathring{\otimes}\left(w_{q+1}^{(c)}+w_{q+1}^{(t)}+w_{q+1}^{(o)}\right)+\left(w_{q+1}^{(c)}+w_{q+1}^{(t)}+w_{q+1}^{(o)}\right) \mathring{\otimes} w_{q+1}\right) . \label{5.5}\tag{4.5}
\end{equation*}
and the commutator error
\begin{equation*}
\mathring{\mathsf{R}}_{\{com\}}:=
\mathcal{R} \mathbb{P}_{H} \operatorname{div}\left(u_{\ell} \mathring{\otimes} u_{\ell}-\left(u_q \mathring{\otimes} u_q\right) *_x \Phi_{\ell} *_t \widetilde{\Phi}_{\ell}\right).
\end{equation*}

Furthermore, it is worthy noting that, as exemplified in references \cite{TBV19,YLZ22}, we have the following mathematical expression:
\begin{equation*}
\mathring{\mathsf{R}}_{q+1}=\mathcal{R} \mathbb{P}_{H} \operatorname{div} \mathring{\mathsf{R}}_{q+1} . \label{5.6}\tag{4.6}
\end{equation*}

\subsection{$C_{t, x}^{1}$-Estimate Of Reynolds Stress.}\ \ \ \
Concerning the Reynolds stress estimates \eqref{2.8} in $C_{t, x}^{1}$-norm. By the identity \eqref{5.6}, Sobolev's embedding $W_x^{1,4} \hookrightarrow L_x^{\infty}$,
$$
\begin{aligned}
\left\|\mathring{\mathsf{R}}_{q+1}\right\|_{C_t C_x^1} & \lesssim\left\|\mathcal{R} \mathbb{P}_H\left(\operatorname{div} \mathring{\mathsf{R}}_{q+1}\right)\right\|_{C_t W_x^{2,4}} \\
& \lesssim\left\|\partial_t u_{q+1}+(-\Delta)^{\alpha} u_{q+1}+\operatorname{div}\left(u_{q+1} \otimes u_{q+1}\right)\right\|_{C_t W_x^{1,4}}.
\end{aligned}
$$
Then, using the interpolation inequality (see \cite{HB18}), and \eqref{347}, we get
\begin{align*}
\left\|\mathring{\mathsf{R}}_{q+1}\right\|_{C_t C_x^1} \lesssim & \left\|u_{q+1}\right\|_{C_{t, x}^2}+\left\|u_{q+1}\right\|_{C_t L_x^4}^{\frac{4-2 \alpha}{5}}\left\|u_{q+1}\right\|_{C_t W_x^{5,4}}^{\frac{2 \alpha+1}{5}} \\
& +\sum_{0 \leq N_1+N_2 \leq 2}\left\|u_{q+1}\right\|_{C_{t, x}^{N_1}}\left\|u_{q+1}\right\|_{C_{t, x}^{N_2}} \\
\lesssim & \lambda_{q+1}^{16}.\label{47}\tag{4.7}
\end{align*}

Similarly,
\begin{align*}
\left\|\partial_t \mathring{\mathsf{R}}_{q+1}\right\|_{C_{t, x}} & \lesssim\left\|\partial_t^2 u_{q+1}+ \partial_t(-\Delta)^{\alpha} u_{q+1}+\operatorname{div} \partial_t\left(u_{q+1} \otimes u_{q+1}\right)\right\|_{C_t L_x^4} \\
& \lesssim\left\|u_{q+1}\right\|_{C_{t, x}^2}+\left\|u_{q+1}\right\|_{C_{t, x}^4}+\left\|u_{q+1} \otimes u_{q+1}\right\|_{C_{t, x}^2}\\
& \lesssim \lambda_{q+1}^{16}. \label{48}\tag{4.8}
\end{align*}

Hence, the verification of the $C_{t,x}^{1}$-estimate \eqref{2.8} for $\mathring{\mathsf{R}}_{q+1}$ have been established.
\subsection{$L_{t, x}^1$-Decay Of Reynolds Stress.}\ \ \ \
In the following section, we focus on verifying $L_{t, x}^{1}$-decay \eqref{2.10} of the Reynolds stress $\mathring{\mathsf{R}}_{q+1}$. As Calder\'{o}n-Zygmund operators are bounded in the space $L_{x}^{\varrho}$, where $1 < \varrho < 2$, we give a specific choice:
\begin{equation*}
\varrho:=\frac{2-12 \varepsilon}{2-13 \varepsilon} \in(1,2), \label{5.12}\tag{4.9}
\end{equation*}
where $\varepsilon$ is defined in \eqref{2.4}. Note that,
\begin{equation*}
(2-12 \varepsilon)\left(1-\frac{1}{\varrho}\right)=\varepsilon \label{5.13}\tag{4.10}
\end{equation*}
and
\begin{equation*}
(r_{\perp}r_{\|})^{\frac{1}{\varrho}-1}=\lambda^{\varepsilon}, \quad (r_{\perp}r_{\|})^{\frac{1}{\varrho}-\frac{1}{2}}=\lambda^{-1+7\varepsilon}. \label{5.14}\tag{4.11}
\end{equation*}

\noindent {\bf (i) Linear Error.} For the acceleration part of the linear error, taking time derivative of \eqref{4.34} and using identity \eqref{4.9}, we obtain
$$
\begin{aligned}
\partial_{t}\Big(w&_{q+1}^{(p)}+w_{q+1}^{(c)}\Big)  =(\lambda r_{\perp}N_{\Lambda})^{-1} \sum_{k} \nabla^{\perp}\left[\partial_{t}\left(a_{(k)} g_{(k)}\right) \boldsymbol{\Psi}_{k}\right]+(\lambda r_{\perp}N_{\Lambda})^{-1} \sum_{k} \nabla^{\perp}\left[a_{(k)} g_{(k)} \partial_{t} \boldsymbol{\Psi}_{k}\right] \\
& =(\lambda r_{\perp}N_{\Lambda})^{-1} \sum_{k} \nabla^{\perp}\left[\partial_{t}\left(a_{(k)} g_{(k)}\right) \boldsymbol{\Psi}_{k}\right]+(\lambda r_{\perp}N_{\Lambda})^{-1} \mu N_{\Lambda}\sum_{k} \nabla^{\perp}\left[a_{(k)} g_{(k)}\left(k_{1} \cdot \nabla\right) \boldsymbol{\Psi}_{k}\right].
\end{aligned}
$$
Note that $\mathcal{R}
\nabla^{\perp}$ is a Calder\'{o}n-Zygmund operator on $\mathbb{T}^2$, we can employ Lemmas \ref{D.1}, \ref{D.3}, and \ref{D.5} to estimate the first term:
$$
\begin{aligned}
(\lambda r_{\perp}N_{\Lambda}&)^{-1}\left\|\sum_{k} \mathcal{R} \nabla^{\perp}\left[\partial_{t}\left(a_{(k)} g_{(k)}\right) \boldsymbol{\Psi}_{k}\right]\right\|_{L^{1}_{t} L_{x}^{\varrho}}\\
 & \lesssim (\lambda r_{\perp})^{-1} \sum_{k}\left\|a_{(k)}\right\|_{C_{x, t}^{1}}\left\|g_{(k)}\right\|_{W^{1,1}}\left\|\boldsymbol{\Psi}_{k}\right\|_{L_{t}^{\infty} L_{x}^{\varrho}} \\
& \lesssim \ell^{-4}(\lambda r_{\perp})^{-1}\sigma \tau^{\frac{1}{2}} r_{\perp}(r_{\perp} r_{\|})^{\frac{1}{\varrho}-\frac{1}{2}}.
\end{aligned}
$$

As for the second term, we first estimate the derivative of $\boldsymbol{\Psi}_k$ in the direction $k_1$ as
$$
\left\|\left(k_{1} \cdot \nabla\right) \boldsymbol{\Psi}_{k}\right\|_{L_{t}^{\infty} L^{\varrho}_{x}} \lesssim \lambda r_{\perp} \frac{r_{\perp}}{r_{\|}}(r_{\perp}r_{\|})^{\frac{1}{\varrho}-\frac{1}{2}} .
$$
where the order of derivative is $\lambda\frac{r_{\perp}}{r_{\|}}$ (rather than $\lambda$ for the full gradient), together with Lemmas \ref{4.3}, \ref{4.5} imply that
$$
\begin{aligned}
(\lambda r_{\perp}N_{\Lambda})^{-1} \mu N_{\Lambda}&\left\|\sum_{k} \mathcal{R} \nabla^{\perp}\left[a_{(k)} g_{(k)}\left(k_{1} \cdot \nabla\right) \boldsymbol{\Psi}_{k}\right]\right\|_{L^{1}_{t}L_{x}^{\varrho}} \\
& \lesssim (\lambda r_{\perp})^{-1} \mu \sum_{k}\left\|a_{(k)}\right\|_{L_{x, t}^{\infty}}\left\|g_{(k)}\right\|_{L^{1}}\left\|\left(k_{1} \cdot \nabla\right) \boldsymbol{\Psi}_{k}\right\|_{L_{t}^{\infty} L_{x}^{\varrho}} \\
& \lesssim\ell^{-1} \frac{r_{\perp}}{r_{\|}} \mu(r_{\perp} r_{\|})^{\frac{1}{\varrho}-\frac{1}{2}}\tau^{\frac{1}{2}} .
\end{aligned}
$$
Due to \eqref{2.4} and \eqref{4.1}, combining both terms we obtain
\begin{align*}
\left\| \mathcal{R}\partial_{t}\left(w_{q+1}^{(p)}+w_{q+1}^{(c)}\right)
\right\|_{L^{1}_{t}L_{x}^{\varrho}} & \lesssim \ell^{-4}\left(\tau^{\frac{1}{2}} \sigma\lambda^{-1}(r_{\perp} r_{\|})^{\frac{1}{\varrho}-\frac{1}{2}}+\frac{r_{\perp}}{r_{\|}}\mu(r_{\perp} r_{\|})^{\frac{1}{\varrho}-\frac{1}{2}}\tau^{-\frac{1}{2}}\right) \\
& \lesssim \ell^{-4} \lambda^{-5\varepsilon}.\label{5.15}\tag{4.12}
\end{align*}

Considering the viscosity term $(-\Delta)^\alpha w_{q+1}$ and using \eqref{4.40}, we have
\begin{align*}
\left\|\mathcal{R}(-\Delta)^{\alpha} w_{q+1}\right\|_{L_{t}^{1} L_{x}^{\varrho}} \lesssim & \left\|\mathcal{R}(-\Delta)^{\alpha} w_{q+1}^{(p)}\right\|_{L_{t}^{1} L_{x}^{\varrho}}+\left\|\mathcal{R}(-\Delta)^{\alpha} w_{q+1}^{(c)}\right\|_{L_{t}^{1} L_{x}^{\varrho}} \\
& +\left\|\mathcal{R}(-\Delta)^{\alpha} w_{q+1}^{(t)}\right\|_{L_{t}^{1} L_{x}^{\varrho}}+\left\|\mathcal{R}(-\Delta)^{\alpha} w_{q+1}^{(o)}\right\|_{L_{t}^{1} L_{x}^{\varrho}}. \label{5.16}\tag{4.13}
\end{align*}
To estimate the right-hand side mentioned above, we employ the interpolation inequality, Lemma \ref{D.6}, along with the condition that $3-2\alpha \geq 20 \varepsilon$ to deduce
\begin{align*}
\left\|\mathcal{R}(-\Delta)^{\alpha} w_{q+1}^{(p)}\right\|_{L_{t}^{1} L_{x}^{\varrho}} & \lesssim\left\||\nabla|^{2 \alpha-1} w_{q+1}^{(p)}\right\|_{L_{t}^{1} L_{x}^{\varrho}} \\
& \lesssim\left\|w_{q+1}^{(p)}\right\|_{L_{t}^{1} L_{x}^{\varrho}}^{\frac{4-2 \alpha}{3}}\left\|w_{q+1}^{(p)}\right\|_{L_{t}^{1} W_{x}^{3, \varrho}}^{\frac{2 \alpha-1}{3}} \\
& \lesssim \ell^{-1} \lambda^{2 \alpha-1} (r_{\perp} r_{\|})^{\frac{1}{\varrho}-\frac{1}{2}} \tau^{-\frac{1}{2}} \lesssim \ell^{-1} \lambda^{-\varepsilon} , \label{5.17}\tag{4.14}
\end{align*}
and
\begin{align*}
& \left\|\mathcal{R}(-\Delta)^{\alpha} w_{q+1}^{(c)}\right\|_{L_{t}^{1} L_{x}^{\varrho}} \lesssim \ell^{-6} \lambda^{2 \alpha-1} (r_{\perp} r_{\|})^{\frac{1}{\varrho}-\frac{1}{2}} \frac{r_{\perp}}{r_{\|}} \tau^{-\frac{1}{2}} \lesssim \ell^{-6} \lambda^{-9 \varepsilon},  \label{5.18}\tag{4.15}\\
& \left\|\mathcal{R}(-\Delta)^{\alpha} w_{q+1}^{(t)}\right\|_{L_{t}^{1} L_{x}^{\varrho}} \lesssim \ell^{-2} \lambda^{2 \alpha-1} \mu^{-1} (r_{\perp} r_{\|})^{\frac{1}{\varrho}-1} \lesssim \ell^{-2} \lambda^{-2\varepsilon}, \label{5.19} \tag{4.16}\\
& \left\|\mathcal{R}(-\Delta)^{\alpha} w_{q+1}^{(o)}\right\|_{L_{t}^{1} L_{x}^{\varrho}} \lesssim \ell^{-13} \sigma^{-1} \lesssim \ell^{-13} \lambda^{-2 \varepsilon}. \label{5.20}\tag{4.17}
\end{align*}

Therefore, combining \eqref{5.17}-\eqref{5.20} with the observation that $\ell^{-13}\ll\lambda^{\varepsilon}$, we derive that
\begin{equation*}
\left\|\mathcal{R}(-\Delta)^{\alpha} w_{q+1}\right\|_{L_{t}^{1} L_{x}^{\varrho}} \lesssim \ell^{-1} \lambda^{-\varepsilon}. \label{5.21}\tag{4.18}
\end{equation*}

The nonlinear part of \eqref{5.3} remains unestimated. By \eqref{2.6} and Lemma \ref{D.6},
\begin{align*}
& \left\|\mathcal{R} \mathbb{P}_{H} \operatorname{div}\left(w_{q+1} \otimes u_{\ell}+u_{\ell} \otimes w_{q+1}\right)\right\|_{L_{t}^{1} L_{x}^{\varrho}} \\
\lesssim & \left\|w_{q+1} \otimes u_{\ell}+u_{\ell} \otimes w_{q+1}\right\|_{L_{t}^{1} L_{x}^{\varrho}} \\
\lesssim & \left\|u_{\ell}\right\|_{C_{t,x}^{1}}\left\|w_{q+1}\right\|_{L_{t}^{1} L_{x}^{\varrho}} \\
\lesssim & \lambda_{q}^{7}\left(\ell^{-1} (r_{\perp} r_{\|})^{\frac{1}{\varrho}-\frac{1}{2}} \tau^{-\frac{1}{2}}+\ell^{-2}\mu^{-1} (r_{\perp}r_{\|})^{\frac{1}{\varrho}-1}+\ell^{-5} \sigma^{-1}\right) \lesssim \ell^{-6} \lambda^{-2 \varepsilon} . \label{5.22}\tag{4.19}
\end{align*}

Based on the analysis presented in \eqref{5.15}, \eqref{5.21}, and \eqref{5.22}, we arrive at
\begin{equation*}
\left\|\mathring{\mathsf{R}}_{\{l i n\}}\right\|_{L_{t}^{1} L_{x}^{\varrho}} \lesssim \ell^{-4} \lambda^{-5\varepsilon}+\ell^{-1} \lambda^{-\varepsilon}+\ell^{-6} \lambda^{-2 \varepsilon} \lesssim \ell^{-6} \lambda^{-\varepsilon} . \label{5.23}\tag{4.20}
\end{equation*}

\noindent {\bf (ii) Oscillation Error.} Now we are prepared to estimate the oscillation error, we  decompose it into three components as:
$$
\mathring{\mathsf{R}}_{\{osc\}}=\mathring{\mathsf{R}}_{\{o s c .1\}}+\mathring{\mathsf{R}}_{\{o s c .2\}}+\mathring{\mathsf{R}}_{\{o s c .3\}},
$$
where the low-high spatial oscillation error
$$
\mathring{\mathsf{R}}_{\{o s c .1\}}:=\sum_{k \in \Lambda} \mathcal{R} \mathbb{P}_{H} \mathbb{P}_{\neq 0}\left(g_{(k)}^{2} \mathbb{P}_{\neq 0}\left(\mathbf{W}_{k} \otimes \mathbf{W}_{k}\right) \nabla\left(a_{(k)}^{2}\right)\right),
$$
the high temporal oscillation error
$$
\mathring{\mathsf{R}}_{\{o s c .2\}}:=-\mu^{-1} \sum_{k \in \Lambda} \mathcal{R} \mathbb{P}_{H} \mathbb{P}_{\neq 0}\left(\partial_{t}\left(a_{(k)}^{2} g_{(k)}\right) \varphi_{\left(k_{1}\right)} ^{2} (\psi_{(k)}^{\prime})^{2} k_{1}\right),
$$
and the low frequency error
$$
\mathring{\mathsf{R}}_{\{o s c .3\}}:=-\sigma^{-1} \sum_{k \in \Lambda} \mathcal{R} \mathbb{P}_{H} \mathbb{P}_{\neq 0}\left(h_{(k)} \fint_{\mathbb{T}^{2}} \mathbf{W}_{k} \otimes \mathbf{W}_{k} \mathrm{d} x \partial_{t} \nabla\left(a_{(k)}^{2}\right)\right).
$$

In terms of the low-high spatial oscillation error $\mathring{\mathsf{R}}_{\text {osc.1 }}$, note that the velocity flows exhibit high oscillations
$$
\mathbb{P}_{\neq 0}\left(\mathbf{W}_{k} \otimes \mathbf{W}_{k}\right)=\mathbb{P}_{\geq\left(\lambda r_{\perp} / 2\right)}\left(\mathbf{W}_{k} \otimes \mathbf{W}_{k}\right),
$$
we employ Lemmas \ref{D.1} and \ref{D.5}, and use Lemma \ref{H.3} specifically with the substitutions $a=
\nabla \left(a_{(k)}^{2}\right)$ and $f=\varphi_{\left(k_{1}\right)}^{2} (\psi_{(k)}^{\prime})^{2}$ to derive
\begin{align*}
\left\|\mathring{\mathsf{R}}_{\{o s c .1\}}\right\|_{L_{t}^{1} L_{x}^{\varrho}} & \lesssim \sum_{k \in \Lambda}\left\|g_{(k)}\right\|_{L_{t}^{2}}^{2}\left\||\nabla|^{-1} \mathbb{P}_{\neq 0}\left(\mathbb{P}_{\geq\left(\lambda r_{\perp} / 2\right)}\left(\mathbf{W}_{k} \otimes \mathbf{W}_{k}\right) \nabla\left(a_{(k)}^{2}\right)\right)\right\|_{C_{t} L_{x}^{\varrho}} \\
& \lesssim \sum_{k \in \Lambda}\left\||\nabla|^{3}\left(a_{(k)}^{2}\right)\right\|_{C_{t, x}} (\lambda r_{\perp})^{-1}\left\|\varphi_{\left(k_{1}\right)} \right\|^{2}_{C_{t} L_{x}^{2\varrho}}\left\|\psi_{(k)}^{\prime}\right\|^{2}_{C_{t} L_{x}^{2\varrho}} \\
& \lesssim \ell^{-13} \lambda^{-1} r_{\perp}^{\frac{1}{\varrho}-2} r_{\|}^{\frac{1}{\varrho}-1}. \label{5.24}\tag{4.21}
\end{align*}

Furthermore, we utilize Lemmas \ref{D.3}, \ref{D.1}, \ref{D.5} and introduce the large temporal oscillation parameter $\mu$ to balance the high temporal oscillation error $\mathring{\mathsf{R}}_{\{osc. 2\}}$:
\begin{align*}
\left\|\mathring{\mathsf{R}}_{\{o s c .2\}}\right\|_{L_{t}^{1} L_{x}^{\varrho}} & \lesssim \mu^{-1} \sum_{k \in \Lambda}\left\|\mathcal{R}\mathbb{P}_{H} \mathbb{P}_{\neq 0}\left(\partial_{t}\left(a_{(k)}^{2} g_{(k)}^{2}\right) \varphi_{\left(k_{1}\right)} ^{2} (\psi_{(k)}^{\prime})^{2} k_{1}\right)\right\|_{L_{t}^{1} L_{x}^{\varrho}} \\
& \lesssim \mu^{-1} \sum_{k \in \Lambda}\left(\left\|\partial_{t}\left(a_{(k)}^{2}\right)\right\|_{C_{t, x}}\left\|g_{(k)}^{2}\right\|_{L_{t}^{1}}+\left\|a_{(k)}\right\|_{C_{t, x}}^{2}\left\|\partial_{t} (g_{(k)}^{2})\right\|_{L_{t}^{1}}\right)\\&\quad
\left\|\varphi_{\left(k_{1}\right)} \right\|_{C_{t} L_{x}^{2 \varrho}}^{2}\left\|\psi_{(k)}^{\prime}\right\|_{L_{x}^{2 \varrho}}^{2} \\
& \lesssim\left(\ell^{-5}+\ell^{-2}\sigma \tau \right) \mu^{-1} (r_{\perp} r_{\|})^{\frac{1}{\varrho}-1} \\
& \lesssim \ell^{-2}  \sigma \tau\mu^{-1} (r_{\perp} r_{\|})^{\frac{1}{\varrho}-1} . \label{5.25}\tag{4.22}
\end{align*}

The estimates of the low-frequency error $\mathring{\mathsf{R}}_{\{osc. 3\}}$ can be estimated through the utilization of \eqref{4.18} and
Lemma \ref{D.5}:
\begin{align*}
\left\|\mathring{\mathsf{R}}_{\{o s c .3\}}\right\|_{L_{t}^{1} L_{x}^{\varrho}} &\lesssim \sigma^{-1} \sum_{k \in \Lambda}\left\|h_{(k)} \partial_{t} \nabla\left(a_{(k)}^{2}\right)\right\|_{L_{t}^{1} L_{x}^{\varrho}}\\
& \lesssim \sigma^{-1} \sum_{k \in \Lambda}\left\|h_{(k)}\right\|_{C_{t}}\left(\left\|a_{(k)}\right\|_{C_{t, x}}\left\|a_{(k)}\right\|_{C_{t, x}^{2}}+\left\|a_{(k)}\right\|_{C_{t, x}^{1}}^{2}\right)\\
& \lesssim \ell^{-8} \sigma^{-1} . \label{5.26}\tag{4.23}
\end{align*}

Therefore, combining \eqref{5.24}- \eqref{5.26} and utilizing expressions \eqref{4.1}, \eqref{5.14}, along with the constraint $0 < \varepsilon < (3-2\alpha) / 20$, we arrive at
\begin{align*}
\left\|\mathring{\mathsf{R}}_{\{osc\}}\right\|_{L_{t}^{1} L_{x}^{\varrho}} & \lesssim \ell^{-13} \lambda^{-1} r_{\perp}^{\frac{1}{\varrho}-3} r_{\|}^{\frac{1}{\varrho}-1}+\ell^{-2}  \sigma \tau \mu^{-1}  (r_{\perp}r_{\|})^{\frac{1}{\varrho}-1}+\ell^{-8} \sigma^{-1} \\
& \lesssim \ell^{-13} \lambda^{-\varepsilon}+\ell^{-2} \lambda^{2 \alpha-3+15 \varepsilon}+\ell^{-8} \lambda^{-2 \varepsilon} \\
& \lesssim \ell^{-13} \lambda^{-\varepsilon} . \label{5.27}\tag{4.24}
\end{align*}

\noindent {\bf (iii) Corrector Error.} Drawing inspiration from \cite{YLZ22}, we introduce $p_{1}, p_{2} \in(1, \infty)$ satisfying the conditions
$$
\frac{1}{p_{1}}=1-\widetilde{\eta}, \quad \frac{1}{p_{1}}=\frac{1}{p_{2}}+\frac{1}{2},
$$
with $\widetilde{\eta} \leq \varepsilon / (2-12 \varepsilon)$. by H\"{o}lder inequality alongside Lemmas \ref{D.6}, \eqref{4.54}, we deduce
\begin{align*}
\left\|\mathring{\mathsf{R}}_{\{cor\}}\right\|_{L_{t}^{1} L_{x}^{p_{1}}} & \lesssim\left\|w_{q+1}^{(p)} \otimes\left(w_{q+1}^{(c)}+w_{q+1}^{(t)}+w_{q+1}^{(o)}\right)-\left(w_{q+1}^{(c)}+w_{q+1}^{(t)}+w_{q+1}^{(o)}\right) \otimes w_{q+1}\right\|_{L_{t}^{1} L_{x}^{p_{1}}} \\
& \lesssim\left\|w_{q+1}^{(c)}+w_{q+1}^{(t)}+w_{q+1}^{(o)}\right\|_{L_{t}^{2} L_{x}^{p_{2}}}\left(\left\|w_{q+1}^{(p)}\right\|_{L_{t, x}^{2}}+\left\|w_{q+1}\right\|_{L_{t, x}^{2}}\right) \\
& \lesssim \delta_{q+1}^{\frac{1}{2}}\left(\ell^{-6} (r_{\perp} r_{\|})^{\frac{1}{p_{2}}-\frac{1}{2}}\frac{r_{\perp}}{r_{\|}}+\ell^{-2} \mu^{-1} (r_{\perp}r_{\|})^{\frac{1}{p_{2}}-1}\tau^{\frac{1}{2}}+\ell^{-7} \sigma^{-1}\right) \\
& \lesssim \ell^{-7} \delta_{q+1}^{\frac{1}{2}}\left(\lambda^{-8 \varepsilon-2 \widetilde{\eta}-12 \widetilde{\eta} \varepsilon}+\lambda^{2 \widetilde{\eta}-2\varepsilon-12 \widetilde{\eta} \varepsilon}+\lambda^{-2 \varepsilon}\right) \lesssim \ell^{-5} \lambda^{-\varepsilon} \label{5.28}\tag{4.25}
\end{align*}
where the last step is derived from the inequality $2 \widetilde{\eta}-2\varepsilon-12 \widetilde{\eta} \varepsilon \leq-\varepsilon $.

Therefore, combining the estimates \eqref{5.23}, \eqref{5.27}, and \eqref{5.28}, along with the observation that $\ell^{-13}\ll\lambda^{\frac{\varepsilon}{2}}$, we arrive at
\begin{align*}
\left\|\mathring{\mathsf{R}}_{q+1}\right\|_{L_{t, x}^{1}} & \leq\left\|\mathring{\mathsf{R}}_{\{l i n\}}\right\|_{L_{t}^{1} L_{x}^{\varrho}}+\left\|\mathring{\mathsf{R}}_{\{osc\}}\right\|_{L_{t}^{1} L_{x}^{\varrho}}+\left\|\mathring{\mathsf{R}}_{\{cor\}}\right\|_{L_{t}^{1} L_{x}^{p_{1}}} +\left\|\mathring{\mathsf{R}}_{\{com\}}\right\|_{L_{t}^{1} L_{x}^{\varrho}}\\
& \lesssim \ell^{-6} \lambda^{-\varepsilon}+\ell^{-13} \lambda^{-\varepsilon}+\ell^{-5} \lambda^{-\varepsilon} +\ell\lambda_{q}^{14}\\
& \leq \delta_{q+2}. \label{5.29}\tag{4.26}
\end{align*}

Hence, the $L_{t,x}^{1}$-estimate \eqref{2.10} of new Reynolds stress is verified.

In conclusion, the iterative estimates \eqref{2.6}-\eqref{2.13} have been verified in the previous sections, we only need to prove the remaining temporal inductive inclusion \eqref{2.14}. By definitions,
$$
\begin{aligned}
& \operatorname{supp}_t w_{q+1} \subseteq \bigcup_{k \in \Lambda} \operatorname{supp}_t a_{(k)} \subseteq \mathscr{N}_{2\ell}\left(\operatorname{supp}_t \mathring{\mathsf{R}}_{q}\right),
\end{aligned}
$$
which yield that
$$
\begin{aligned}
& \operatorname{supp}_t u_{q+1} \subseteq \operatorname{supp}_t u_{\ell} \cup \operatorname{supp}_t w_{q+1} \subseteq \mathscr{N}_{2\ell}\left(\operatorname{supp}_t u_q \cup \operatorname{supp}_t \mathring{\mathsf{R}}_{q} \right)
\end{aligned}
$$

Moreover, by \eqref{5.2},
$$
\begin{aligned}
& \operatorname{supp}_t\mathring{\mathsf{R}}_{q+1} \subseteq \bigcup_{k \in \Lambda} \operatorname{supp}_t a_{(k)} \cup \operatorname{supp}_t u_\ell,
\end{aligned}
$$
In view of \eqref{212}, we arrive at
$$
\operatorname{supp}_t \mathring{\mathsf{R}}_{q+1} \subseteq \mathscr{N}_{2\ell}\left(\operatorname{supp}_t u_q \cup \operatorname{supp}_t\mathring{\mathsf{R}}_{q}\right) .
$$

Since $2\ell \ll \delta_{q+2}^{\frac{1}{2}}$ due to \eqref{2.3}, the inductive inclusion \eqref{2.14} is verified. Therefore, the proof of Lemma \ref{B.3} is complete.

{\centering
\section{Proof Of Main Theorem}\label{F}}
Based on the establishment of the iteration framework in Lemma \ref{B.3}, we begin to proof our main theorem:\\
{\bf Proof Of Theorem \ref{A.2}.}  Taking $u_{0}=\tilde{u}$ and set
\begin{align*}
\mathring{\mathsf{R}}_{0} & :=\mathcal{R}\left(\partial_{t} u_{0}+(-\Delta)^{\alpha} u_{0}\right)+u_{0} \mathring{\otimes} u_{0},  \tag{5.1}\\
P_{0} & :=-\frac{1}{2}\left|u_{0}\right|^{2} , \tag{5.2}
\end{align*}
then the pair $\left(u_{0}, \mathring{\mathsf{R}}_{0}\right)$ constitutes one solution of the initial step to the equations \eqref{2.1}. We introduce the notation $\delta_{1}=\|\mathring{\mathsf{R}}_{0}\|_{L_{t,x}^{1}}$ and select the parameter $a$ to be sufficiently large, ensuring that the conditions \eqref{2.6}-\eqref{2.10} hold true at the base level $q=0$. Subsequently, by Lemma \ref{B.3}, we can construct a sequence of solutions $(u_{q}, \mathring{\mathsf{R}}_{q})_{q}$ to \eqref{2.1} that satisfy the inductive estimates \eqref{2.6}-\eqref{2.14} for all $q\geq0$.

By using the interpolation, \eqref{2.2}, \eqref{2.6} and \eqref{2.11}, we infer that for any $\beta^{\prime} \in\left(0, \frac{\beta}{7+\beta}\right)$,
\begin{align*}
\sum_{q \geq 0}\left\|u_{q+1}-u_{q}\right\|_{H_{t, x}^{\beta^{\prime}}} & \leq \sum_{q \geq 0}\left\|u_{q+1}-u_{q}\right\|_{L_{t, x}^{2}}^{1-\beta^{\prime}}\left\|u_{q+1}-u_{q}\right\|_{H_{t, x}^{1}}^{\beta^{\prime}} \\
& \lesssim \sum_{q \geq 0} \delta_{q+1}^{\frac{1-\beta^{\prime}}{2}} \lambda_{q+1}^{7 \beta^{\prime}} \\
& \lesssim \delta_{1}^{\frac{1-\beta^{\prime}}{2}} \lambda_{1}^{7 \beta^{\prime}}+\sum_{q \geq 1} \lambda_{q+1}^{-\beta\left(1-\beta^{\prime}\right)+7 \beta^{\prime}}<\infty. \tag{5.3}
\end{align*}

By virtue of \eqref{4.64}, we have
\begin{equation*}
\sum_{q \geq 0}\left\|u_{q+1}-u_{q}\right\|_{L_{t}^{\gamma} L_{x}^{p}}<\infty. \tag{5.4}
\end{equation*}

Therefore, $\left\{u_{q}\right\}_{q \geq 0}$ is established as a Cauchy sequence within the super-critical space $H_{t, x}^{\beta^{\prime}}\cap L_{t}^{\gamma} L_{x}^{p}$, and there exists $u$ such that $\lim _{q \rightarrow \infty}\left(u_{q}\right)=u$ in $H_{t, x}^{\beta^{\prime}}\cap L_{t}^{\gamma} L_{x}^{p}$. Additionally, considering the fact that $\lim _{q \rightarrow \infty} \mathring{\mathsf{R}}_{q}=0$ in $L_{t, x}^{1}$, we conclude that $u\in H_{t, x}^{\beta^{\prime}}\cap L_{t}^{\gamma} L_{x}^{p}$ is a weak solution to \eqref{1.1}. This verifies the assertions $(i)$ and $(ii)$.

Concerning the $\varepsilon_{*}$-neighborhood of temporal supports, combining the temporal inductive inclusion \eqref{2.14}
and taking into account $\sum_{q \geq 0} \delta_{q+2}^{1 / 2} \leq \varepsilon_{*}$ for $a$ large enough, then the temporal support statement $(iii)$ is verified. Finally, for the $\varepsilon_{*}$-close between the solution $u$ and the given vector field $\tilde{u}$ in $L_{t}^{1} L_{x}^{2}\cap L_{t}^{\gamma} W_{x}^{s, p}$, by \eqref{2.12} and \eqref{2.13}, $(iv)$ is easy to be obtained.

Therefore, the proof of Theorem \ref{A.2} is complete.\\

{\centering
\section{Appendix}\label{G}}

In this section, we shall exhibit some preliminary results that have been previously utilized in the preceding sections.

\begin{Lemma}\label{G.1} (\cite{ACL22}, Geometric Lemma 3.1). There exists a set $\Lambda \subset \mathbb{S} \cap \mathbb{Q}^{2}$ consisting vectors $k$, and positive smooth functions $\gamma_{(k)}: B_{C_{\mathsf{R}}}(\mathrm{Id}) \rightarrow \mathbb{R}$, such that for $\mathsf{R} \in B_{C_{\mathsf{R}}}(Id)$, it holds that
\begin{equation*}
\mathsf{R}=\sum_{k \in \Lambda} \gamma_{(k)}^{2}(\mathsf{R}) k_{1} \otimes k_{1},\tag{6.1}
\end{equation*}
where $(k, k_{1})$ is an orthonormal basis and $B_{C_{\mathsf{R}}}(\mathrm{Id})$ represents the ball with a radius of $C_{\mathsf{R}}>0$ centered on the identity matrix in the $2 \times 2$ symmetric matrix space.
\end{Lemma}

As stated in the previous work \cite{BCV22}, there exists  $N_{\Lambda}\in \mathbb{N}$ that fulfills
\begin{equation*}
\left\{N_{\Lambda} k, N_{\Lambda} k_{1}\right\} \subseteq N_{\Lambda} \mathbb{S} \cap \mathbb{Z}^{2}. \label{7.2}\tag{6.2}
\end{equation*}
We introduce the notation $M_{*}$ to denote a geometric constant that satisfies the following inequality:
\begin{equation*}
\sum_{k \in \Lambda}\left\|\gamma_{(k)}\right\|_{C^{4}\left(B_{C_{\mathsf{R}}}(\mathrm{Id})\right)} \leq M_{*}. \tag{6.3}
\end{equation*}

Subsequently, we revisit the $L^{p}$ decorrelation lemma, introduced by Lemma 2.4 in \cite{AC21} (see also \cite{BV19}, Lemma 3.7). This lemma serves as a crucial component in obtaining the $L_{t, x}^{2}$ estimates for the perturbations.
\begin{Lemma}\label{H.2}(\cite{AC21}, Lemma 2.4). For any $p \in[1, \infty]$, there exist $\sigma \in \mathbb{N}$ and smooth functions $f, g: \mathbb{T}^{d} \rightarrow \mathbb{R}$ such that,
\begin{equation*}
\left|\|f g(\sigma \cdot)\|_{L^{p}\left(\mathbb{T}^{d}\right)}-\|f\|_{L^{p}\left(\mathbb{T}^{d}\right)}\|g\|_{L^{p}\left(\mathbb{T}^{d}\right)}\right| \lesssim \sigma^{-\frac{1}{p}}\|f\|_{C^{1}\left(\mathbb{T}^{d}\right)}\|g\|_{L^{p}\left(\mathbb{T}^{d}\right)}. \tag{6.4}
\end{equation*}
\end{Lemma}

The following stationary phase lemma, as a primary tool, plays a pivotal role in addressing the errors associated with Reynolds stress.
\begin{Lemma}\label{H.3} (\cite{LTP20}, Lemma 7.4). For any given $1<p<\infty$, $\lambda\in \mathbb{Z}_{+}$, $a \in C^{2}\left(\mathbb{T}^{2}, \mathbb{R}\right)$ and $f\in L^{p}\left(\mathbb{T}^{2}, \mathbb{R}^{2}\right)$, one has
$$
\left\||\nabla|^{-1} \mathbb{P}_{\neq 0}\left(a \mathbb{P}_{\geq \lambda} f\right)\right\|_{L^{p}\left(\mathbb{T}^{2}\right)} \lesssim \lambda^{-1}\left\|\nabla^{2} a\right\|_{L^{\infty}\left(\mathbb{T}^{2}\right)}\|f\|_{L^{p}\left(\mathbb{T}^{2}\right)}
.$$
\end{Lemma}

\section*{Acknowledgment.}
Zhong Tan thanks the support by National Natural Science Foundation of China Grant (12071391, 12231016) and the Guangdong Basic and Applied Basic Research Foundation (2022A1515010860). Xinliang Li is supported by National Natural Science Foundation of
China Grant (12401297) and the project funded by China Postdoctoral Science Foundation (2024T170577, 2023M742401).
\section*{Conflicts of interest}

The authors declare that they have no conflict of interest.

\section*{Data availability}

Data sharing is not applicable to this article as no data were created or analyzed in this study.

\begin{center}

\end{center}
\end{document}